\documentclass[reqno]{amsart}
\usepackage{amsfonts,amsmath,amssymb,color}

\numberwithin{equation}{section}

\usepackage[babel=true,kerning=true]{microtype}

\usepackage{amsthm,leqno}
\usepackage{hyperref}

\usepackage{amsfonts,amsmath,amssymb,enumerate}
\usepackage{graphicx}
\usepackage{tikz}
\usetikzlibrary{decorations.markings}
\usetikzlibrary{plotmarks}
\usetikzlibrary{patterns}

\newtheorem{theo}{Theorem}
\newtheorem{cl}{Claim}
\newtheorem*{lem}{Lemma}
\newtheorem*{cor}{Corollary}

\newcommand{\eps}{\epsilon}

\begin{document}

\title[Maximal metrics for the first eigenvalue]{Existence and regularity of maximal metrics for the first Laplace eigenvalue on surfaces}
\author{Romain Petrides} 
\address{Romain Petrides, Universit\'e de Lyon, CNRS UMR 5208, Universit\'e Lyon 1, Institut Camille Jordan, 43 bd du 11 novembre 1918, F-69622 Villeurbanne cedex, France.}
\email{romain.petrides@univ-lyon1.fr}

\begin{abstract} We investigate in this paper the existence of a metric which maximizes the first eigenvalue of the Laplacian on Riemannian surfaces. We first prove that, in a given conformal class, there always exists such a maximizing metric which is smooth except at a finite set of conical singularities. This result is similar to the beautiful result concerning Steklov eigenvalues recently obtained by Fraser and Schoen \cite{FRA2013}. Then we get existence results among all metrics on surfaces of a given genus, leading to the existence of minimal  isometric immersions of smooth compact Riemannian manifold $(M,g)$ of dimension $2$ into some $k$-sphere by first eigenfunctions. At last, we also answer a conjecture of Friedlander and Nadirashvili \cite{FRIED1999} which asserts that the supremum of the first eigenvalue of the Laplacian on a conformal class can be taken as close as we want of its value on the sphere on any orientable surface. 
\end{abstract}

\maketitle

Let $\left(\Sigma,g\right)$ be a smooth compact Riemannian surface without boundary. The eigenvalues of the Laplacian $\Delta_g = - div_g\left(\nabla\right)$ form a discrete sequence 
$$0=\lambda_0<\lambda_1\left(\Sigma,g\right)\le \lambda_2\left(\Sigma,g\right)\le \dots$$
Getting bounds on these eigenvalues depending on the metric or the topology of $\Sigma$ has been the subject of intensive studies in the past decades. In this paper, we shall focus on the first eigenvalue $\lambda_1$. One can for instance consider the first conformal eigenvalue of $\left(\Sigma,g\right)$ defined by 
\begin{equation}\label{defLambda1conf}
\Lambda_1\left(\Sigma,[g]\right)= \sup_{\tilde{g}\in [g]} \lambda_1\left(\tilde{g}\right) Vol_{\tilde{g}}\left(\Sigma\right)\hskip.1cm.
\end{equation}
If one looks at the infimum of the first eigenvalue in a given conformal class, it is always $0$. Now one can also study invariants which depend only on the topology of the surface. For orientable surfaces, one can define for any genus $\gamma\ge 0$
\begin{equation}\label{defLambda1}
\Lambda_1\left(\gamma\right) = \sup_{g} \lambda_1(g)Vol_g\left(\Sigma\right) = \sup_{[g]} \Lambda_1\left(\Sigma,[g]\right)
\end{equation}
where $\Sigma$ is a compact orientable surface of genus $\gamma$. One can also look at 
$$\inf_{[g]}  \Lambda_1\left(\Sigma,[g]\right)\hskip.1cm.$$

\medskip Natural questions about these quantities are to get explicit values or explicit bounds on it, and whether or not the supremum (or infimum) in their definition is achieved by some metric and, if yes, how regular these extremal metrics are. Yang and Yau \cite{YAN1980} (see also \cite{LIY1982}) obtained an upper-bound for the first eigenvalue of the Laplacian on a surface, depending only on the genus $\gamma$ of the surface. In case of orientable surfaces, this reads as
\begin{equation} \label{yangyau} 
\Lambda_1\left(\gamma\right) \leq 8\pi \left[\frac{\gamma+3}{2}\right]\hskip.1cm. 
\end{equation}
Colbois and El Soufi \cite{COL2003} gave an explicit lower bound of $\Lambda_1\left(\Sigma,[g]\right)$ on 
any closed Riemannian surface and proved that 
$$\Lambda_1\left(\Sigma,[g]\right) \ge \Lambda_1\left(\mathbb{S}^2,[\hbox{can}]\right)$$ 
and by the work of Hersch \cite{HER1970}, we know that $\Lambda_1\left(\mathbb{S}^2,[\hbox{can}]\right) = 8\pi$. A lower bound for $\Lambda_1\left(\gamma\right)$ can be obtained from \cite{BROOKS2001} and \cite{BBD1988} (see \cite{FRS2013})~:
\begin{equation} \label{lowerbound} 
\Lambda_1\left(\gamma\right) \geq \frac{3\pi}{4}\left(\gamma-1\right)\hskip.1cm. 
\end{equation}

\medskip Exact values of these quantities were obtained for small genus and for specific conformal classes. Let us mention the sphere (Hersch  \cite{HER1970}), the projective plane (Li-Yau \cite{LIY1982}), the torus (Girouard \cite{GIR2009bis} and Nadirashvili \cite{NAD1996}), the Klein bottle (El Soufi-Giacomini-Jazar \cite{ELS2006} and Jakobson-Nadirashvili-Polterovitch \cite{JAK2006}), the genus $2$ surfaces (Jakobson-Levitin-Nadirashvili-Nigam-Polterovitch \cite{JAK2005}, Karpukhin \cite{KAR2013}).

\medskip Concerning $\Lambda_1\left(\Sigma,[g]\right)$, we prove the following theorem~:

\begin{theo}\label{thm1}  Let $(\Sigma,g)$ be a compact Riemannian surface without boundary. Then
$$ \Lambda_1\left(\Sigma,[g]\right) > \Lambda_1\left(\mathbb{S}^2,[can]\right)= 8\pi$$
if $\Sigma$ is not diffeomorphic to $\mathbb{S}^2$. Moreover, there is an extremal 
metric $\tilde{g} \in [g]$, smooth except maybe at a finite number of points corresponding to conical
singularities, such that $\Lambda_1\left(\Sigma,[g]\right) = \lambda_1\left(\tilde{g}\right)Vol_{\tilde{g}}\left(\Sigma\right)$.
\end{theo}

This theorem contains a rigidity result which states that the sphere is characterized by having the minimal first conformal eigenvalue. It also contains an existence result of "smooth" maximal metrics. Note that, on the sphere, we know since the work of Hersch \cite{HER1970} that maximal metrics exist and are all smooth since they consist in all metrics isometric to the standard one. As observed in \cite{KOK2011}, conical singularities naturally appear for extremal metrics. Indeed, the conformal factor relating $\tilde{g}$ to $g$ is $\left|\nabla \Phi\right|_g^2$ where $\Phi$ is some smooth harmonic map from $M$ into some sphere $\mathbb{S}^k$. The zeros of $\left|\nabla \Phi\right|_g^2$ are isolated and correspond to branch points of the harmonic map $\Phi$ as proved in Salamon \cite{SAL1985}. The case of genus $2$ surface treated in \cite{KAR2013} show that extremal metrics indeed possess sometimes conical singularities. In this respect, our existence result seems completely optimal. In \cite{KOK2011}, Kokarev proved that any maximizing sequence of metrics, provided that our rigidity result was true, converges to a Radon measure without atoms. He then got some partial regularity results on this measure.  Note that, here, we do not prove that any maximizing sequence converges to a "smooth" maximizer, which may not be true. We select carefully a maximizing sequence which converges to a "smooth" maximizer. In \cite{NADSIRE2010}, assuming that ${\displaystyle \Lambda_1\left(\Sigma,[g]\right) > 8\pi}$, which is by now a consequence of our result, the authors announced the existence of a maximizer with a rather different proof we do not fully understand. 

\medskip Note also that, by Kokarev \cite{KOK2010}, and thanks to the rigidity part of our theorem, we know that the set of "smooth" maximizers given by our theorem is compact as soon as $M$ is not diffeomorphic to the sphere $\mathbb{S}^2$. On the sphere $\mathbb{S}^2$, this compactness result is of course false.

\medskip Capitalizing on this first existence result, we are also able to obtain the following~:

\begin{theo}\label{thm2}
Let $\Sigma$ be a compact orientable Riemannian surface without boundary of genus $\gamma\geq 1$. If $ \Lambda_1\left(\gamma\right) > \Lambda_1\left(\gamma-1\right)$, then $\Lambda_1(\gamma)$ is achieved by a metric which is smooth except at a finite set of conical singularities.
\end{theo}

Note that the case of the sphere is already treated in Hersch \cite{HER1970}. Note also that the case of the torus ($\gamma=1$) is already known~: we have $\Lambda_1\left(1\right)=\frac{8\pi^2}{\sqrt{3}}$ and the maximal metric is given by the flat equilateral torus (see \cite{NAD1996}). At last, the genus $2$ case was recently obtained~: we have $\Lambda_1\left(2\right)=16\pi$ and there is a family of maximal metrics (see \cite{KAR2013}). 

\medskip The spectral gap $\Lambda_1\left(\gamma\right)>\Lambda_1\left(\gamma-1\right)$ necessarily holds for an infinite number of $\gamma$ thanks to the lower bound (\ref{lowerbound}). It is believed to hold for all genuses. The extremal metric in the theorem is the pull-back of the induced metric of a minimal immersion (with branched points) of $\Sigma$ into some sphere ${\mathbb S}^k$. As a classical corollary of the above theorem, we obtain the following~:

\begin{cor}\label{corollary}
If $\gamma\ge 1$ and if $\Lambda_1\left(\gamma\right)>\Lambda_1\left(\gamma-1\right)$, which is the case at least for an infinite number of $\gamma$, there exists a minimal immersion (possibly with branch points) of a compact surface $\Sigma$ of genus $\gamma$ into some sphere ${\mathbb S}^k$ by first eigenfunctions.
\end{cor}

There have been lot of works about minimal immersions of surfaces into spheres. In particular, they are necessarily given by eigenfunctions (not only first eigenfunctions) thanks to Takahashi \cite{TAK1966}. For existence results of such immersions, we refer to two classical papers by Lawson \cite{LAW1970} and Bryant \cite{BRY1982}. Concerning minimal embeddings in ${\mathbb S}^3$, it is conjectured by Yau \cite{YAU1982} that they all come from first eigenfunctions (see \cite{BRE2013} and \cite{CHO2006} for recent surveys on this subject). However, minimal immersions by first eigenfunctions are not so numerous.  For instance, it has been proved by Montiel and Ros \cite{MON1986} that there is at most one minimal immersion by first eigenfunctions in any given conformal class. In the case of genus $1$, it was also proved by El Soufi and Ilias \cite{ELS2000} that the only minimal immersions by first eigenfunctions of the torus are the Clifford torus (in ${\mathbb S}^3$) and the flat equilateral torus (in ${\mathbb S}^5$). So our corollary is interesting because it provides an infinite number of new minimal immersions into spheres by first eigenfunctions.

\medskip At last, we prove a conjecture stated in \cite{FRIED1999} about the infimum of the first conformal eigenvalue on any orientable surface~:

\begin{theo}\label{thm3}
Let $\Sigma$ be a smooth compact orientable surface. Then 
$$\inf_{[g]}\Lambda_1\left(\Sigma,[g]\right)=8\pi$$
and this infimum is never attained except on the sphere.
\end{theo}

This result had already been proved in \cite{GIR2009bis} in genus $1$ but was left open in higher genuses up to now. 

\medskip The paper is organized as follows :

\medskip We first prove in section \ref{rigidity} the rigidity part of theorem \ref{thm1}. The idea of the proof goes back to Ledoux \cite{LED2000} and Druet \cite{DRU2002} in higher dimensions. We start from some Moser-Trudinger type  inequality (see \cite{CHE1979, DIN1997, MOS1971}) which possesses extremal functions. These extremal functions are excellent candidates to provide conformal factors for which the new metric has a large $\lambda_1$. However, we have to deal with some degeneracy problems which could occur. 

\medskip Then, we prove the existence of a "smooth" extremal metric for $\Lambda_1\left(M,[g]\right)$. In section \ref{nonconcentration}, we prove some fine non-concentration estimates for sequences of unit volume metrics in a given conformal class with large first eigenvalue. This non-concentration phenomenon was first observed by Girouard \cite{GIR2009bis} and Kokarev \cite{KOK2011}. Section \ref{construction} is devoted to the construction of our specific maximizing sequence, following ideas of Fraser and Schoen \cite{FRA2013} when dealing with the Steklov eigenvalue problem. It is obtained by solving a regularized maximization problem. We derive a fine Euler-Lagrange characterization for this new variational problem. This leads to a maximizing sequence of smooth metrics for which the first eigenspace possesses nice properties. Section \ref{convergence} makes an intensive use of the non-concentration estimates of section \ref{nonconcentration} to get finer and finer estimates on these first eigenfunctions. This permits then to pass to the limit and to prove theorem \ref{thm1}. 

\medskip Section \ref{sectionminimalimmersion} is devoted to the proof of theorem \ref{thm2}. Since we already have the existence of a maximizing metric in any given conformal class thanks to theorem \ref{thm1}, it remains to prove that the supremum among all conformal classes is achieved. For that purpose, we pick up a sequence of maximizing conformal classes and prove that this sequence does not degenerate. We follow ideas of Zhu \cite{ZHU2010} who made a careful study of sequence of harmonic maps into spheres on hyperbolic surfaces which degenerate.

\medskip The last section is devoted to the proof of theorem \ref{thm3}. It is in some sense similar to the proof of theorem \ref{thm2}, except that we just have to construct a sequence of fully degenerating conformal classes $\left(c_\alpha\right)$ on some hyperbolic surface $\Sigma$ for which we prove that $\Lambda_1\left(\Sigma,c_\alpha\right)\to 8\pi$.

\medskip\textbf{Acknowledgements}

It is a pleasure to thank my thesis advisor O. Druet for his support and comments during the preparation of this work, and P. Laurain for discussions about the asymptotics of harmonic maps.

\section{The rigidity result}\label{rigidity}
\subsection{Extremal functions for a sharp Moser-Trudinger inequality}

\medskip Let $(M,g)$ be a smooth compact Riemannian surface with $Vol_g(M) = 1$, not diffeomorphic to $\mathbb{S}^2$. The aim of this section is to find a metric $\tilde{g}\in [g]$ such that $\lambda_1\left(\tilde{g}\right)>8\pi$. We consider the functional $J_g$ defined for $u\in H_1^2\left(M\right)$ by 
\begin{equation}\label{defJ}
J_g(u)= \frac{1}{4\pi} \int_M \left|\nabla u\right|^2_g dv_g + 2\int_M u dv_g - \log\left( \int_M e^{2u} dv_g \right)\hskip.1cm.
\end{equation}
Cherrier \cite{CHE1979} proved that 
\begin{equation}\label{defalphag}
-\alpha(g)=\inf_{u\in H_1^2\left(M\right)} J_g(u)
\end{equation}
is always finite. Note that $J_g$ is translation invariant with respect to the constant functions. 

\medskip In the following, up to a harmless conformal change of the metric, we assume that $g$ is the metric in $[g]$ with volume 1 and constant Gaussian curvature $K_g\equiv K_0$. Since $M$ is not diffeomorphic to the sphere, we know that $K_0\le 2\pi$. Then there exists a ball $B_\delta(0)\subset C^2\left(M\right)$ centered at $0$ and of radius $\delta>0$ such that for any $v\in B_\delta(0)$, 
$$K_{e^{2v}g}<4\pi\hskip.1cm.$$
We claim that there exists $v_0\in B_\delta(0)\subset C^2\left(M\right)$ such that $8\pi$ is not in the spectrum of $\Delta_{e^{2v_0}g}$. It is completely obvious since the spectrum is discrete and it scales like the volume to the power $-1$. Thus there exists an open set $V\subset C^2\left(M\right)$ such that, for any $v\in V$,
\begin{equation}\label{cond1}
K_{e^{2v}g}<4\pi
\end{equation}
and 
\begin{equation}\label{cond2}
8\pi\hbox{ is not in the spectrum of }\Delta_{e^{2v}g}\hskip.1cm.
\end{equation}
Let us now define the following map~:
$$\begin{array}{ccccc}
F & : & W^{2,2}(M) \times V & \longrightarrow & L^2(M) \\
& & (u,v) & \longmapsto & e^{-2v}\Delta_{g} u + 4\pi - 4\pi e^{2u}  \\
\end{array}$$
It is well defined because of the Sobolev embedding $W^{2,2}(M) \subset \mathcal{C}^0\left(M\right)$ and $F$ is a $\mathcal{C}^{\infty}$ map.

\medskip By a result of \cite{DIN1997} and thanks to (\ref{cond1}), for any $v\in V$, the functional $J_{e^{2v}g}$ attains its minimum. Let $u\in H_1^2\left(M\right)$ be such that $J_{e^{2v}g}(u)=-\alpha\left(e^{2v}g\right)$ normalized by ${\displaystyle \int_M e^{2u}e^{2v}\, dv_g=1}$. Then $u\in C^3(M)$ and satisfies the Euler-Lagrange equation 
\begin{equation}\label{DJ}
F(u,v)= e^{-2v} \Delta _g u + 4\pi - 4\pi e^{2u}=0\hskip.1cm.
\end{equation}
Moreover, we have that  
$$
D^2 J_{e^{2v}g}(u)\left(\varphi,\varphi\right) = \frac{1}{2\pi} \left(\int_M \left\vert \nabla \varphi\right\vert_{\tilde{g}}^2\, dv_{\tilde{g}} - 8\pi \int_M \varphi^2\, dv_{\tilde{g}} + 8\pi \left(\int_M \varphi\, dv_{\tilde{g}}\right)^2\right)
\ge 0
$$
where $\tilde{g}=e^{2(u+v)} g$. This means in particular that $\lambda_1\left(e^{2(u+v)}g\right)\ge 8\pi$. Since 
\begin{equation}\label{eqFred}
D_u F(u,v)\left(\varphi\right) = e^{2u}\left(\Delta_{e^{2(u+v)}g} \varphi-8\pi \varphi\right)\hskip.1cm,
\end{equation}
we have $\lambda_1\left(e^{2(u+v)}g\right)> 8\pi$ as soon as $D_u F(u,v)\left(\varphi\right)$ is invertible. Thus, in order to prove the rigidity part of the theorem, we just need to find $v\in V$ such that $D_u F(u,v)$ is invertible for all solutions $u$ of $F(u,v)=0$. 

\medskip For that purpose, we shall apply the following theorem of Fredholm theory (see for instance \cite{HEN2005}, Theorem 5.4, page 63) to our function $F$~:

\begin{theo} \label{gene} Let $X$,$Y$ two separable Banach spaces, $U$ an open set of $X$, $V$ a separable $\mathcal{C}^{\infty}$ Banach manifold and $F\in\mathcal{C}^{\infty}(U\times V,Y)$ which satisfies :
\begin{itemize}
\item For all $(u,v)\in F^{-1}(0)$, $DF(u,v)$ is surjective.
\item For all $(u,v)\in F^{-1}(0)$, $D_u F(u,v)$ is a Fredholm operator.
\end{itemize}
Then, there exists a countable intersection of open dense sets (a residual set) $\Sigma \subset V$ such that for all $v\in\Sigma$, and for all $u\in F(.,v)^{-1}(0)$, $D_u F(u,v)$ is surjective. 
\end{theo}

\medskip By (\ref{eqFred}), it is clear that $D_u F(u,v)$ is a Fredholm operator. It remains to prove that if $(u,v)\in F^{-1}(0)$, $DF(u,v)$ is surjective. We have 
$$ DF(u,v).(\theta,\tau) = e^{2u} \left( \Delta_{\tilde{g}} \theta - 8\pi  \theta\right) - 2 \tau e^{-2v} \Delta_g u$$
where $\tilde{g}=e^{2(u+v)}g$. Since the image of $DF(u,v)$ contains a finite codimensional closed space, it is a closed space. Moreover, the $L^2$ norms induced by the metrics $g$ and $e^{2v}g$ are equivalent. Then, it suffices to prove that the orthogonal of the image is $0$ in $L^2(e^{2v}g)$. Assume on the contrary that there exists $\phi\in L^2(M)$, $\phi\not\equiv 0$, such that
$$ \forall \theta \in W^{2,2}(M), \int_M \left(\Delta_{\tilde{g}} \theta - 8\pi \theta \right)\phi dv_{\tilde{g}} =0 $$
and
$$ \forall \tau \in \mathcal{C}^2(M),  \int_M \phi \tau \Delta_g u dv_g = 0 \hskip.1cm.$$
The first condition implies that $\phi$ is an eigenfunction for $\tilde{g}$ with eigenvalue $8\pi$. We deduce from the second condition knowing that $e^{-2v} \Delta_g u + 4\pi = 4\pi e^{2u}$ that 
$$ \phi \left( e^{2u} - 1 \right) = 0 \hskip.1cm.$$
Since $\phi$ is non-zero on a dense set of points by the maximum principle, $e^{2u} \equiv 1$. This implies that $\tilde{g} = e^{2v}g$ has a $8\pi$ eigenvalue, which contradicts (\ref{cond2}). 

\medskip Thus we can apply the above theorem to our function $F$ to get the existence of some $v\in V$ such that $D_u F(u,v)$ is surjective for all $u\in W^{2,2}\left(M\right)$ such that $F(u,v)=0$. As already said, this ends the proof of the rigidity part of theorem \ref{thm1}. For such a $v$ and for a minimal function $u$ for $J_{e^{2v}g}$, we get that $\lambda_1\left(e^{2(u+v)}g\right)>8\pi$.

\section{Existence of maximal metrics in a conformal class}

Let $(M,g)$ be a smooth compact unit volume Riemannian surface without boundary. We choose for all the proof some $\delta>0$, some $C_0>1$, a family $\left(x_i\right)_{i=1,\dots, N}$ of points in $M$ and smooth functions $v_i:M\mapsto {\mathbb R}$ such that 
\begin{itemize}
\item for any $i\in \left\{1,\dots,N\right\}$, the Gauss curvature of $g_i=e^{2v_i}g$ is $0$ in the ball $B_{g_i}\left(x_i, 2\delta\right)=\Omega_i$ so that, in the exponential chart for the metric $g_i$ at $x_i$, the metric $g_i$ is the Euclidean metric. 

\item ${\displaystyle M=\bigcup_{i=1}^N \omega_i}$ where ${\displaystyle \omega_i = B_{g_i}\left(x_i, \delta\right)}$.

\item for any $i\in \left\{1,\dots,N\right\}$, $C_0^{-2} \le e^{2v_i}\le C_0^2$ in $\Omega_i$.

\end{itemize}

\medskip\noindent Note that, for any $i\in \left\{1,\dots,N\right\}$,
\begin{equation}\label{eq-partition1}
B_{g}\left(x,C_0^{-1}r\right) \subset B_{g_i}\left(x,r\right)\subset B_g\left(x, C_0r\right)\hbox{ for all }x\in \omega_i\hbox{ and all } 0<r\le \delta\hskip.1cm.
\end{equation}
During all this section, in order to get uniform estimates, we may assume, without loss of generality that every sequence $\{x_{\epsilon}\}$ of points of $M$ lies in some $\omega_i$ where $i$ is fixed. Indeed, every subsequence of $\{x_{\epsilon}\}$ has a subsequence which satisfies this property.

\subsection{Non-concentration estimates of metrics with high first eigenvalue}\label{nonconcentration}

In this subsection, we let $\{ e^{2u_{\epsilon}}g \}_{\epsilon}$ be a sequence of unit volume metrics such that 
\begin{equation}\label{higheigenvalue}
\lambda_1\left(e^{2u_{\epsilon}}g\right) \geq 8\pi + \alpha\hbox{ for some }\alpha>0 \hbox{ fixed.}
\end{equation}
Note that we also know that ${\displaystyle \lambda_1\left(e^{2u_{\epsilon}}g\right)\le \Lambda_1\left(M,[g]\right)}$. By Kokarev \cite{KOK2011}, lemma 2.1 and lemma 3.1, the following non-concentration result follows from (\ref{higheigenvalue})~:

\begin{cl}[Kokarev \cite{KOK2011}]\label{claim-Kokarev}
Assume that $\left(u_{\epsilon}\right)$ is a sequence of smooth functions on $M$ such that (\ref{higheigenvalue}) holds. Then 
$$\lim_{r\to 0} \limsup_{\eps\to 0} \sup_{x\in M} \int_{B_{g}\left(x,r\right)} e^{2u_\epsilon}\, dv_g =0\hskip.1cm.$$
\end{cl}

We set for $\Omega$ an open subset of $M$ and $\mu\in\mathcal{M}(M)$, the set of Radon measures on $M$,  
$$ \lambda_{*}(\Omega,\mu) = \inf_{\phi\in\mathcal{C}^{\infty}_c(\Omega)} \frac{\int_{\Omega} \left|\nabla\phi\right|^2_g \, dv_g}{\int_{\Omega} \phi^2 \, d\mu} $$ 
and
$$ \beta(\Omega,\mu) = \sup\left\{ \frac{\mu(F)}{Cap_2\left(F,\Omega\right)} ; F \subset \Omega \text{ is a compact set} \right\} $$
where, for a compact set $F\subset \Omega$, $Cap_2\left(F,\Omega\right)$ is the variational capacity of $(F,\Omega)$ defined by 
$$ Cap_2(K,\Omega) = \inf \left\{ \int_{\Omega} \left|\nabla \phi\right|_g^2\, dv_g ; \phi \in \mathcal{C}^{\infty}_c(\Omega) , \phi = 1 \text{ on K} \right\}\hskip.1cm. $$
Isocapacitary inequalities proved in \cite{MAZ1985}, section 2.3.3, corollary of theorem 2.3.2, give that 
\begin{equation} \label{isoc} 
\frac{1}{4\beta(\Omega,\mu)} \leq \lambda_{*}(\Omega,\mu) \leq \frac{1}{\beta(\Omega,\mu)}\hskip.1cm. \end{equation}
Thanks to these capacity estimates, we can refine the non-concentration result of Kokarev and obtain a quantitative one~:

\begin{cl} \label{cap} 
Assume that $\left(u_{\epsilon}\right)$ is a sequence of smooth functions on $M$ such that (\ref{higheigenvalue}) holds. Then there exists $C_1>0$ such that 
$$\int_{B_g\left(x,r\right)} e^{2u_{\epsilon}}dv_g \leq \frac{C_1}{\log{\frac{1}{r}}} $$
for all $\epsilon>0$ and all $r>0$. 
\end{cl}

\medskip {\bf Proof.} We first prove that there exists $r_0 >0$ such that for any $0<r\le r_0$, 
\begin{equation} \label{conc} 
\forall \epsilon >0, \forall x\in M, \frac{1}{\lambda_{*}(B_g\left(x,r_0\right),e^{2u_{\epsilon}}g)} \leq \frac{2}{\lambda_{1}(M,e^{2u_{\epsilon}}g)} \leq \frac{1}{4\pi}\hskip.1cm.
\end{equation}
Indeed, choose $\psi_{\epsilon} \in E_1\left(B_g\left(x,r_0\right),e^{2u_{\epsilon}}g\right)$ with ${\displaystyle \int_{M} \psi_{\epsilon}^2 e^{2u_{\epsilon}}dv_g =1}$ and let us write that 
\begin{eqnarray*}
&&\int_{B_{g}\left(x,r_0\right)}\psi_{\epsilon}^2e^{2u_{\epsilon}}dv_g - \left(\int_{B_g\left(x,r_0\right)}\psi_{\epsilon}e^{2u_{\epsilon}}dv_g\right)^2 \\
&&\quad \leq \frac{1}{\lambda_1\left(M,e^{2u_{\epsilon}}g\right)} \int_{M}\left|\nabla\psi_{\epsilon}\right|_g^2\, dv_g = \frac{\lambda_{*}\left(B_g\left(x,r_0\right),e^{2u_{\epsilon}}g\right)}{\lambda_1\left(M,e^{2u_{\epsilon}}g\right)}\hskip.1cm. 
\end{eqnarray*}
By H\"older's inequality, we deduce that 
$$\int_{B_g\left(x,r_0\right)} \psi_{\epsilon}^2 e^{2u_{\epsilon}}dv_g \left(1-\int_{B_g\left(x,r_0\right)}e^{2u_{\epsilon}}dv_g\right) \leq \frac{\lambda_{*}(B_g\left(x,r_0\right),e^{2u_{\epsilon}}g)}{\lambda_1(M,e^{2u_{\epsilon}}g)}\hskip.1cm.$$
By claim \ref{claim-Kokarev}, there exists $r_0 >0$ such that
$$ \forall \epsilon>0, \forall x \in M, \int_{B_g\left(x,r_0\right)}e^{2u_{\epsilon}}dv_g \leq \frac{1}{2} $$
and (\ref{conc}) follows for this $r_0$. 

\medskip Let's fix 
$$r_1 = \frac{1}{C_0^2} \min\left\{\delta,r_0\right\}\hskip.1cm.$$
Then, for any $i\in \left\{1,\dots,N\right\}$, any $x\in \omega_i$, by (\ref{eq-partition1}) and (\ref{conc}), we have that 
\begin{eqnarray*}
\lambda_\star\left(B_{g_i}\left(x,r_1\right), e^{2u_\eps}g_i\right) &\ge & \lambda_\star\left(B_{g}\left(x,C_0 r_1\right), e^{2u_\eps}g_i\right)\\
&\ge & C_0^{-2}\lambda_\star\left(B_{g}\left(x,C_0r_1\right), e^{2u_\eps}g\right)\\
&\ge & \frac{4\pi}{C_0^2}\hskip.1cm.
\end{eqnarray*}
Writing thanks to (\ref{isoc}) that 
$$\int_{B_{g_i}\left(x,r\right)} e^{2u_\epsilon}\, dv_g \le \frac{Cap_2\left(B_{g_i}\left(x,r\right),B_{g_i}\left(x,r_1\right)\right)}{\lambda_\star\left(B_{g_i}\left(x,r_1\right), e^{2u_\eps}g_i\right)}$$
for all $0<r<r_1$ and thanks to the fact that $g_i$ is isometric to the Euclidean metric that
$$Cap_2\left(B_{g_i}\left(x,r\right),B_{g_i}\left(x,r_1\right)\right) = \frac{2\pi}{\ln\frac{r_1}{r}}\hskip.1cm,$$
we get that 
$$\int_{B_{g_i}\left(x,r\right)} e^{2u_\epsilon}\, dv_g \le \frac{C_0^2}{2\ln\frac{r_1}{r}}$$
for all $0<r<r_1$. This clearly leads to the conclusion of the claim. \hfill $\diamondsuit$

\medskip We now focus on the eigenfunctions associated to the first non-zero eigenvalue of such a sequence of metrics. We will prove that the nodal sets of such eigenfunctions can not concentrate to a point~:

\begin{cl}\label{nod}
There exists $\delta_1>0$ such that for any $\eps>0$, any $f\in E_1\left(e^{2u_\eps}g\right)$ and any $x\in M$,
$$f\left(x\right)=0 \Rightarrow \exists\, y\in \partial B_g\left(x,\delta_1\right)\hbox{ s.t. }f(y)=0\hskip.1cm.$$
\end{cl}

\medskip {\bf Proof.} Assume by contradiction that for any $j\in\mathbb{N}$, there exist $\epsilon_j>0$, $x_j\in M$ and $f_j\in E_1\left(e^{2u_{\epsilon_j}} g\right)$ such that 
\begin{equation} \label{zer} 
f_{j}\left(x_j\right)=0 \hbox{ and } \forall y\in \partial B_g\left(x_j,2^{-j}\right), f_j(y)\neq 0\hskip.1cm. \end{equation}
Then, by the maximum principle, $f_j$ changes sign in $B_g\left(x_j,2^{-j}\right)$. By the Courant nodal theorem (see \cite{COU1953}), there are two connected nodal domains $D_j^{1}$ and $D_j^{2}$ for $f_j$. We know that for $m\in\{1,2\}$, $D_j^m \cap B_{2^{-j}}(x_j) \neq \emptyset$. Therefore, thanks to (\ref{zer}) and what we just said, there is one of this nodal domain, let's say $D_j^1$, which satisfies $D_j^1\subset B_{2^{-j}}(x_j)$.

\medskip Then $f_j$ is an eigenfunction on $(D_j^2,e^{2u_{\epsilon_j}})$ with Dirichlet boundary conditions. Since $f_j>0$, it is an eigenfunction for $\lambda_{*}\left(D_j^2,e^{2u_{\epsilon_j}}\right)=\lambda_{\eps_j}$. Up to a subsequence, $x_j \to x\in M$ as $j\to\infty$. Thanks to claim \ref{claim-Kokarev}, there is an open neighborhood of $\{x\}$, $B$, with 
$$ \int_{B} e^{2u_{\epsilon_j}} \leq \frac{1}{2}\hskip.1cm. $$
Since $Cap_2\left(\{x\},B\right)=0$, we can find $\psi \in \mathcal{C}^{\infty}_c(M\setminus \{x\})$ such that $\psi=1$ on $M\setminus B$ and
$$ \int_M \left|\nabla \psi\right|^2_gdv_g \leq 2\pi\hskip.1cm.$$
We write then that 
$$ \lambda_{\epsilon_j} = \lambda_{*}\left(D_j^2,e^{2u_{\epsilon_j}}\right) \leq \frac{\int_{D_j^2}\left|\nabla \psi\right|^2_gdv_g}{\int_{D_j^2} \psi^2 e^{2u_{\epsilon_j}}dv_g} \leq 4\pi $$
which contradicts (\ref{higheigenvalue}). This ends the proof of the claim. \hfill $\diamondsuit$

\medskip At last, the non-concentration estimates just obtained gives the $W^{1,2}(M)$-boundedness of a sequence of normalized eigenfunctions~:

\begin{cl} \label{soboun} 
Any sequence of eigenfunctions $f_{\epsilon}\in E_1\left(e^{2u_{\epsilon}}g\right)$ such that $\int_{M} f_{\epsilon}^2e^{2u_{\epsilon}}dv_g = 1$ is bounded in $W^{1,2}(M)$.
\end{cl}

\medskip {\bf Proof.} We already know that ${\displaystyle \int_M \left|\nabla f_{\epsilon}\right|_g^2dv_g = \lambda_1\left(e^{2u_{\epsilon}}g\right)}$ is bounded. We now prove that $\{e^{2u_{\epsilon}}dv_g\}_{\epsilon}$ is a bounded sequence in $W^{-1,2}(M) = W^{1,2}(M)^{*}$.

\medskip Let us consider a finite covering of $M$ by balls $B_g\left(y_j,r_0\right)$, $j=1,\dots, L$ where $r_0>0$ is given by (\ref{conc}) and let $\left(\psi_j\right)$ be a partition of unity associated to this covering. 
For $\psi \in W^{1,2}(M)$, we have that
\begin{eqnarray*} 
\int_M \psi e^{2u_{\epsilon}} dv_g &=& \sum_{j=1}^{L} \int_{B_g\left(y_j,r_0\right)} \psi\psi_j e^{2u_{\epsilon}} dv_g \\
&\leq& \sum_{j=1}^{L} \left(\int_{B_g\left(y_j,r_0\right)} (\psi_j \psi)^2 e^{2u_{\epsilon}}dv_g \right)^{\frac{1}{2}} \left(\int_{B_g\left(y_j,r_0\right)}e^{2u_{\epsilon}}dv_g \right)^{\frac{1}{2}} \\
&\le & \sum_{j=1}^L \frac{1}{\lambda_{*}\left(B_g\left(y_j,r_0\right),e^{2u_{\epsilon}}g\right) ^{\frac{1}{2}}} \left(\int_{M} \left|\nabla(\psi_j\psi)\right|^2 dv_g \right)^{\frac{1}{2}}\\
&\le&  D_0\left\|\psi\right\|_{W^{1,2}(M)}
\end{eqnarray*}
where $D_0$ is independent of $\psi$ and $\epsilon$. Then $\{e^{2u_{\epsilon}}dv_g\}_{\epsilon}$ is a bounded sequence in $W^{-1,2}(M)$ and we get the following Poincar\'e inequality (see \cite{ZIE1989}, lemma 4.1.3)~: there exists $D_2>0$ such that 
$$ \forall \epsilon>0, \forall \psi\in\mathcal{C}^{\infty}(M), \int_{M}\left(\psi - \int_M \psi e^{2u_{\epsilon}}dv_g\right)^2 dv_g \leq D_2 \int_M \left|\nabla\psi\right|^2_g dv_g\hskip.1cm. $$
We apply this inequality to $\psi = f_{\epsilon}$ which has zero mean value with respect to $e^{2u_{\epsilon}}dv_g$ since $f_{\epsilon}$ is a first eigenfunction with respect to $e^{2u_{\epsilon}}g$. We get that
$$\int_{M}f_{\epsilon}^2 dv_g \leq D_2 \int_M \left|\nabla f_{\epsilon}\right|^2_g dv_g=D_2\lambda_1\left(e^{2u_{\epsilon}}g\right)$$
which gives the desired conclusion. \hfill $\diamondsuit$

\subsection{Construction of a maximizing sequence}\label{construction}
\medskip For $\epsilon>0$ and $x,y\in M$, we denote by $p_{\epsilon}(x,y)$ the heat kernel of $(M,g)$ at time $\epsilon$. We let $\mathcal{M}(M)$ be the set of positive Radon measures provided with the weak* topology and $\mathcal{M}_1(M)$ be the subset of probability measures. For $\nu \in\mathcal{M}(M)$, $f\in L^1\left(M,g\right)$ and $\epsilon >0$, we set
$$ K_{\epsilon}\left[\nu\right](x) = \int_M p_{\epsilon}(x,y)\, d\nu(y) $$ 
and
$$ K_{\epsilon}\left[f\right](x) = \int_M p_{\epsilon}(x,y) f(y)\, dv_g(y) $$
so that 
$$\int_M K_\epsilon\left[f\right](x)\, d\nu(x) = \int_M f(x)K_\eps\left[\nu\right]\, dv_g(x)\hskip.1cm.$$
Let us now define the maximizing sequence we will consider. For $\epsilon>0$, we set 
\begin{equation}\label{maxprob}
\lambda_{\epsilon} = \sup_{\nu\in\mathcal{M}_1(M)} \lambda_1\left(K_{\epsilon}[\nu]g\right)\hskip.1cm. 
\end{equation}
Since $K_\epsilon[\nu]>0$ and $K_\eps[\nu]\in C^\infty$, $\lambda_{\epsilon} \leq \Lambda_1\left(M,[g]\right)$. Moreover, since $K_{\epsilon}:\mathcal{M}_1(M) \mapsto \mathcal{C}^k(M)$ is continuous for all $k\geq 0$, every maximizing sequence for $\lambda_{\epsilon}$ converges in $\mathcal{M}_1(M)$. So let $\nu_{\epsilon} \in \mathcal{M}_1(M)$ be such that 
$$ \lambda_{\epsilon} = \lambda_1\left(K_{\epsilon}[\nu_{\epsilon}] g\right)\hskip.1cm.$$
We claim that 
\begin{equation}\label{limlambdaepsilon}
\lim_{\epsilon\to 0} \lambda_\eps = \Lambda_1\left(M,[g]\right)\hskip.1cm.
\end{equation}
We already know that $\lambda_{\epsilon} \leq \Lambda_1\left(M,[g]\right)$ for all $\eps>0$. Let $\eta>0$ and pick up $\tilde{g}\in [g]$ such that $Vol_{\tilde{g}}(M)=1$ and ${\displaystyle \lambda_1\left(\tilde{g}\right) > \Lambda_1\left(M,[g]\right) - \frac{\eta}{2}}$. Since $\lambda_1\left(K_\eps\left[dv_{\tilde{g}}\right]g\right)\to \lambda_1\left(\tilde{g}\right)$ as $\eps\to 0$, there exists $\eps_0>0$ such that
$$\lambda_\eps\ge \lambda_1\left(K_\eps\left[dv_{\tilde{g}}\right]g\right)\ge \lambda_1\left(\tilde{g}\right)-\frac{\eta}{2}\ge  \Lambda_1\left(M,[g]\right)-\eta$$
for all $0<\eps<\eps_0$. This proves (\ref{limlambdaepsilon}).

\medskip We let in the following
$$K_{\epsilon}\left[\nu_{\epsilon}\right]=e^{2u_{\epsilon}}$$
and we have that 
$$\lambda_\eps= \lambda_1\left(e^{2u_\eps}g\right)\to  \Lambda_1\left(M,[g]\right)\hbox{ as }\eps\to 0\hskip.1cm.$$
Let us exploit the fact that $\nu_\eps$ solves the maximization problem (\ref{maxprob})~:

\begin{cl} \label{eul} Let's fix $\epsilon>0$. Then there exist $\left(\phi^1_{\epsilon},\cdots,\phi^{k(\epsilon)}_{\epsilon}\right) \in \mathcal{C}^{\infty}(M,\mathbb{R}^{k})$ such that 
\begin{itemize}
\item $\forall i \in \{1,\cdots,k(\eps)\}$, ${\displaystyle \Delta_g \phi^i_{\epsilon} = \lambda_{\epsilon} e^{2u_{\epsilon}} \phi^i_{\epsilon}}$.

\item  ${\displaystyle \int_M e^{2u_{\epsilon}} \left|\Phi_{\epsilon}\right|^2 \, dv_g = 1}$.

\item ${\displaystyle K_{\epsilon}\left[\left|\Phi_{\epsilon}\right|^2\right] \geq 1}$ on $M$.

\item ${\displaystyle K_{\epsilon} \left[\left|\Phi_{\epsilon}\right|^2\right] = 1}$ on $\text{supp}\left(\nu_{\epsilon}\right)$.

\end{itemize}
Here $\Phi_\eps = \left(\phi^1_{\epsilon},\cdots,\phi^{k(\epsilon)}_{\epsilon}\right)$ and ${\displaystyle \left\vert \Phi_\eps\right\vert^2 = \sum_{i=1}^{k(\eps)} \left(\phi_\eps^{i}\right)^2}$. 
\end{cl}

\medskip {\bf Proof.} Since $\epsilon$ is fixed, up to the end of the proof of this claim, we omit the $\epsilon$ indices of $\lambda_{\epsilon}$, $\nu_{\epsilon}$ and $e^{2u_{\epsilon}}$. 

Let $\mu \in \mathcal{M}(M)$ and $t>0$. We let 
$$\lambda_{t} = \lambda_1\left(K_{\epsilon}[\nu+t\mu] g\right)\hskip.1cm.$$
Note that $\lambda = \lambda_{t=0}$. By continuity, $\lambda_{t} \to \lambda$ as $t \to 0^{+}$. We first prove that 
\begin{equation} \label{der} \lim_{t \to 0^{+}} \frac{\lambda_{t}-\lambda}{t} = \inf_{\phi\in E_1(e^{2u}g)} \left( -\lambda\frac{\int_M K_{\epsilon}\left[\phi^2\right] d\mu}{\int_M \phi^2 e^{2u} dv_g} \right) \hskip.1cm.\end{equation}
We let $\phi \in E_1\left(K_{\epsilon}[\nu] g\right) = E_1\left(e^{2u}g\right)$ and we write that 
\begin{eqnarray*}
\lambda_t \left(\int_M \left(\phi - \frac{\int_M \phi K_\eps[\nu+t\mu]\, dv_g}{\int_M K_\eps[\nu+t\mu]\, dv_g}\right)^2 K_\eps[\nu+t\mu]\, dv_g\right) &\le& \int_M \left\vert \nabla \phi\right\vert_g^2\, dv_g\\
&=& \lambda\int_M e^{2u}\phi^2\, dv_g\hskip.1cm.
\end{eqnarray*}
Since ${\displaystyle K_\eps[\nu+t\mu]= e^{2u}+ t K_\eps[\mu]}$, we easily get that 
$$\lambda_t \left(\int_M \phi^2 e^{2u}\, dv_g + t \int_M \phi^2 K_\eps[\mu]\, dv_g + o(t)\right) \le \lambda\int_M e^{2u}\phi^2\, dv_g$$
so that 
\begin{eqnarray*}
\frac{\lambda_t-\lambda}{t} &\le& -\lambda\frac{\int_M \phi^2 K_\eps[\mu]\, dv_g}{\int_M \phi^2 e^{2u}\, dv_g}+o(1)\\
&=& -\lambda \frac{\int_M K_\eps\left[\phi^2\right] \, d\mu}{\int_M \phi^2 e^{2u}\, dv_g}+o(1)\hskip.1cm.
\end{eqnarray*}
So far we have proved that 
\begin{equation}\label{der1}
\limsup_{t\to 0^+} \frac{\lambda_t-\lambda}{t} \le \inf_{\phi\in E_1(e^{2u}g)} \left( -\lambda\frac{\int_M K_{\epsilon}\left[\phi^2\right] d\mu}{\int_M \phi^2 e^{2u} dv_g} \right)\hskip.1cm.
\end{equation}
Let now $\phi_t \in E_1\left(K_{\epsilon}[\nu+t\mu] g\right)$ with $\left\|\phi_t\right\|_{L^2(K_{\epsilon}[\nu+t\mu] g)}=1$. We have that 
\begin{equation} \label{eqe} 
\Delta_g \phi_t = \lambda_t K_{\epsilon}[\nu+t\mu] \phi_t =\lambda_t \left(e^{2u}+t K_\eps[\mu]\right)\phi_t\hskip.1cm.
\end{equation}
For all $t>0$, $L^2(K_{\epsilon}[\nu+t\mu] g)$ and $L^2(e^{2u}g)$ are the same sets and define equivalent norms and the constants in the equivalence are independent of $t$, by the regularity properties of the heat equation. Then, up to the extraction of a subsequence, by standard elliptic regularity, there exists $\phi \in E_1\left(e^{2u} g\right)$ such that $\phi_t \to \phi$ in $\mathcal{C}^m$ as $t \to 0^{+}$ and $\left\|\phi\right\|_{L^2(e^{2u} g)}=1$. We denote by $\Pi$ the orthogonal projection on $E_1\left(e^{2u} g\right)$ with respect to the $L^2(e^{2u} g)$-norm. Then we rewrite (\ref{eqe}) as 
\begin{equation} \label{eqd} \Delta_g\left(\frac{\phi_t - \Pi \phi_t}{\alpha_t}\right) - \lambda e^{2u} \frac{\phi_t - \Pi \phi_t}{\alpha_t} = \frac{\lambda_{t}-\lambda}{\alpha_t} e^{2u} \phi_t + \frac{t}{\alpha_t} \lambda_{t} K_{\epsilon}[\mu]\phi_t \end{equation}
where 
$$\alpha_t = \left\| \phi_t - \Pi \phi_t \right\|_{\infty} + t + (\lambda - \lambda_t)\hskip.1cm.$$
 Then, up to the extraction of a subsequence we have that 
$$t_0 = \lim_{t\to 0^{+}} \frac{t}{\alpha_t} \hbox{ and } \delta_0 = \lim_{t\to 0^{+}} \frac{\lambda_t - \lambda}{\alpha_t}\hskip.1cm.$$
By standard elliptic theory, up to the extraction of a subsequence, 
\begin{equation} \label{reort} 
\frac{\phi_t - \Pi \phi_t}{\alpha_t}\to R_0 \hbox{ in }C^2\left(M\right)\hbox{ as }t\to 0^+\end{equation}
where $R_0 \in E_1(e^{2u} g)^{\perp} $. Passing to the limit in equation (\ref{eqd}), we get that \begin{equation} \label{eqdl} \Delta_g R_0 - \lambda e^{2u} R_0 = \delta_{0} e^{2u} \phi + t_0 \lambda K_{\epsilon}[\mu] \phi \end{equation}
and 
\begin{equation} \label{eqdln} \left\|R_0\right\|_{\infty} + t_0+ \delta_0 = 1 \hskip.1cm.\end{equation} 
Testing (\ref{eqdl}) againt $\phi$ and using the fact that $R_0 \in E_1(e^{2u} g)^{\perp}$ give that 
$$\delta_0 \int_M e^{2u}\phi^2\, dv_g = - t_0\lambda \int_M K_\eps[\mu]\phi^2\, dv_g =  - t_0\lambda \int_M K_\eps\left[\phi^2\right]\, d\mu\hskip.1cm.$$
If $t_0=0$, then $\delta_0=0$ and then $R_0\equiv 0$ thanks to (\ref{eqdl}) and the fact that $R_0 \in E_1(e^{2u} g)^{\perp}$. This is absurd thanks to (\ref{eqdln}). Thus $t_0\neq 0$ and 
$$\lim_{t\to 0^+} \frac{\lambda_t-\lambda}{t}=\frac{\delta_0}{t_0}= - \lambda \frac{\int_M K_\eps\left[\phi^2\right]\, d\mu}{ \int_M e^{2u}\phi^2\, dv_g }\hskip.1cm.$$
This together with (\ref{der1}) gives (\ref{der}). 

\medskip Now, with a renormalization, $(1+t\int_M d\mu) \lambda_t \leq \lambda$ for all $t\geq 0$ and we deduce from (\ref{der}) that

\begin{equation} \label{ineq} \forall \mu \in \mathcal{M}(M), \exists\, \phi \in E_1(e^{2u}g)\hbox{ s.t. } \int_M \phi^2 e^{2u}\, dv_g =1\hbox{ and } \int_M \left(1-K_{\epsilon}[\phi^2]\right)d\mu \leq 0 \hskip.1cm.\end{equation}
Let us define the following subsets of $\mathcal{C}^{0}(M)$~:
$$K = \left\{ \psi \in \mathcal{C}^0(M) ; \exists \,\phi_1,\cdots,\phi_k \in E_1\left(e^{2u}g\right), \psi = \sum_{i=1}^{k} K_{\epsilon}\left[\phi_i^2\right] - 1 , \int_M \psi \, d\nu = 0 \right\}$$
and
$$F = \left\{ f\in \mathcal{C}^0(M) ; f \geq 0 \right\}\hskip.1cm. $$
The set $F$ is closed and convex. The set $K$ is clearly convex since it is the translation of the convex hull of 
$$C=\left\{K_\eps\left[\phi^2\right]; \phi\in E_1\left(e^{2u}g\right), \, \left\Vert \phi\right\Vert_{L^2\left(e^{2u}g\right)}=1\right\}\hskip.1cm.$$
Since $E_1$ is finite-dimensional, the vector space spanned by $C$ is finite-dimensional and $C$ is compact. Caratheodory's theorem gives that $K$ is also compact. 

If $F\cap K=\emptyset$, by Hahn-Banach theorem, there exists $\mu\in  \mathcal{M}(M)$ such that 
\begin{equation}\label{hb1} \forall f\in F, \int_M f d\mu \geq 0\end{equation}
and 
\begin{equation}\label{hb2} \forall \psi\in K, \int_M \psi d\mu < 0\hskip.1cm. \end{equation}
Then, $\mu$ is a non-zero (by (\ref{hb2})) positive (by (\ref{hb1})) measure and for this measure, (\ref{hb2}) contradicts (\ref{ineq}).

Thus $F\cap K \neq \emptyset$ and there exists $\phi^1,\cdots,\phi^k \in E_1(e^{2u}g)$ with
$$ \int_M \left|\Phi\right|^2 e^{2u}dv_g =1 \hbox{ and } K_{\epsilon}\left[\left|\Phi\right|^2\right] \geq 1 $$
where $\Phi =\left(\phi^1,\dots,\phi^k\right)$. Moreover, we can write that 
$$ 1 =\int_M\left|\Phi\right|^2 e^{2u}dv_g  = \int_M K_{\epsilon}\left[\left|\Phi\right|^2\right] d\nu \geq \int_M d\nu = 1\hskip.1cm.$$
Therefore, $K_{\epsilon}\left[\left|\Phi\right|^2\right] = 1$ $\nu$-a.e. and since $K_{\epsilon}\left[\left|\Phi\right|^2\right]$ is continuous, $K_{\epsilon}\left[\left|\Phi\right|^2\right]=1$ on supp($\nu$). This ends the proof of the claim. \hfill $\diamondsuit$

\medskip Note that the number of eigenfunctions $k(\epsilon)$ depends on $\epsilon$ but since the multiplicity of eigenvalues is bounded by a constant which only depends on the topology of the surface (see \cite{CHE1976}), up to the extraction of a subsequence, we assume in the following that $k$ is fixed.

\subsection{Estimates on eigenfunctions}\label{convergence}

We assume now that $(M,g)$ is not diffeomorphic to the sphere. Then, by the rigidity result proved in section \ref{rigidity}, $\Lambda_1\left(M,[g]\right) > 8\pi$ and we can use the non-concentration estimates of the specific maximizing sequence $\left\{e^{2u_{\epsilon}} g\right\}$ of the previous sections. We denote by $\nu$ the weak$^{*}$ limit of $\{ e^{2u_{\epsilon}} dv_g\}_{\epsilon > 0}$. Then $\nu$ is also the weak$^{*}$ limit of $\{d\nu_{\epsilon}\}$. Indeed, for $\zeta \in \mathcal{C}^{0}(M)$,
\begin{eqnarray*}
\left| \int_M \zeta \left(d\nu_{\epsilon} - e^{2u_{\epsilon}} dv_g \right) \right| &=& \left| \int_M \zeta \left(d\nu_{\epsilon} - K_{\epsilon}[\nu_{\epsilon}] dv_g \right) \right| \\
&=& \left| \int_M \left(\zeta - K_{\epsilon}[\zeta]\right) d\nu_{\epsilon} \right| \leq \left\Vert \zeta - K_{\epsilon}[\zeta]\right\Vert_{\infty} \underset{\epsilon \to 0}{\longrightarrow} 0\hskip.1cm. 
\end{eqnarray*}
We aim at proving that $\nu$ is absolutely continuous with respect to $dv_g$ with a smooth density.

Up the end of the proof, we get finer and finer estimates on the sequence of eigenfunctions given by claim \ref{eul}. For that purpose, we shall use the uniform estimates of the heat kernel $p_{\epsilon}$ on $M$ as $\epsilon \to 0$ (see \cite{BER1971}) :
\begin{equation} \label{heat} p_{\epsilon}(x,y) = \frac{1}{4\pi\epsilon} e^{-\frac{d_g(x,y)^2}{4\epsilon}} + O\left(e^{-\frac{1}{\sqrt{\epsilon}}}\right)\hskip.1cm. \end{equation}

\begin{cl} \label{bou} For all $i\in\{1,\cdots,k\}$, the sequence $\{\phi^i_{\epsilon}\}_{\epsilon}$ is uniformly bounded. 
\end{cl}

\medskip {\bf Proof.} Let $i\in\{1,\cdots,k\}$. Let $\left(x_\eps\right)$ be a sequence of points in $M$ such that $\phi_\eps^i\left(x_\eps\right)=\sup_{M}\left\vert \phi_\eps^i\right\vert$. Up to change $\phi_\eps^i$ into $-\phi_\eps^i$, such a $x_\eps$ does exist. We set 
$$\delta_\eps = d_g\left(x_\eps,\text{supp}\left(\nu_\eps\right)\right)\hskip.1cm.$$
Up to the extraction of a subsequence, we can assume that $x_\eps\in \omega_l$ for some $l$ fixed. We set then 
$$\tilde{\phi}_\eps = \phi_\eps^i\left(\exp_{g_l,x_l}\left(x\right)\right)$$
so that 
$$\Delta_\xi \tilde{\phi}_\eps = \lambda_\eps e^{2\tilde{u}_\eps} \tilde{\phi}_\eps \hbox{ in }B_0\left(2\delta\right)$$
with 
$$\tilde{u}_\eps = \left(u_\eps+ v_l\right)\left(\exp_{g_l,x_l}\left(x\right)\right)\hskip.1cm.$$
We also let in the following 
$$\tilde{x}_\eps = \exp_{g_l,x_l}^{-1}\left(x_\eps\right)\hskip.1cm.$$

\medskip We divide the proof into three cases.

\medskip {\sc Case 1} - We assume that $\delta_\eps^{-1}=O(1)$.

\medskip Then, by the convergence properties of the heat kernel, $\left\{ e^{2u_{\epsilon}}\right\}$ is uniformly bounded in $B_{g}\left(x_{\epsilon},\frac{\delta_{\epsilon}}{2}\right)$ (see (\ref{heat})) and by the claim \ref{soboun}, $\left\{\phi_{\epsilon}^i \right\}$ is bounded in $L^2(M)$. By standard elliptic theory, $\left\{\phi_{\epsilon}^i(x_{\epsilon})\right\}$ is bounded. 

\medskip {\sc Case 2} - We assume that $\delta_\eps = O\left(\sqrt{\epsilon}\right)$. 

\medskip We let 
$$\hat{\phi}_\eps(x) = \tilde{\phi}_\eps\left(\sqrt{\eps}x+\tilde{x}_\eps\right)$$
for $x\in B_0\left(\delta\eps^{-\frac{1}{2}}\right)$. Then 
$$\Delta_\xi \hat{\phi}_\eps =\eps \lambda_\eps   e^{2\tilde{u}_\eps\left(\sqrt{\eps}x+\tilde{x}_\eps\right)} \hat{\phi}_\eps$$
in $B_0\left(\delta\eps^{-\frac{1}{2}}\right)$. By (\ref{heat}), $\left(\epsilon p_{\epsilon}\right)$ is uniformly bounded so that $\left(\eps e^{2\tilde{u}_\eps\left(\sqrt{\eps}x+\tilde{x}_\eps\right)}\right)$ is uniformly bounded. 

Now, thanks to the assumption made in this case and to claim \ref{eul}, there exists $y_\eps\in M$ such that $d_g\left(x_\eps,y_\eps\right)=O\left(\sqrt{\eps}\right)$ and such that ${\displaystyle K_{\epsilon}\left[\left|\Phi_{\epsilon}\right|^2\right](y_{\epsilon}) = 1}$. Let us write then thanks to (\ref{heat}) that 
\begin{eqnarray*}
1= K_{\epsilon}\left[\left|\Phi_{\epsilon}\right|^2\right]\left(y_{\epsilon}\right)&\ge &
K_{\epsilon}\left[\left|\phi^i_{\epsilon}\right|^2\right]\left(y_{\epsilon}\right)\\
&=& \int_M p_\eps\left(x,y_\eps\right) \left(\phi_\eps^i(y)\right)^2\, dv_g(y)\\
&\ge & \frac{1}{4\pi \eps} e^{-R^2}\int_{B\left(y_\eps,2R\sqrt{\eps}\right)}  \left(\phi_\eps^i(y)\right)^2\, dv_g(y) \\
&&\quad+ O\left(e^{-\frac{1}{\sqrt{\eps}}}\right)\int_{B\left(y_\eps,2R\sqrt{\eps}\right)}  \left(\phi_\eps^i(y)\right)^2\, dv_g(y) \hskip.1cm.
\end{eqnarray*}
We let $\tilde{y}_\eps = \frac{1}{\sqrt{\eps}}\left(\exp_{g_l,x_l}^{-1} \left(y_\eps\right) - \tilde{x}_{\eps}\right)$ so that, up to a subsequence, $\tilde{y}_\eps\to y_0$ as $\eps\to 0$ and we deduce from the previous inequality that, for any $R>0$, there exists $C_R>0$ such that
$$\int_{\mathbb{D}_R\left(y_0\right)} \hat{\phi}_\eps^2\, dx \le C_R\hskip.1cm.$$
Thus, by standard elliptic theory, it is clear that $\left\{\hat{\phi}_\eps\right\}$ is uniformly bounded in any compact subset of ${\mathbb R}^2$. This gives in this second case that $\left\{\phi_\eps^i\left(x_\eps\right)\right\}$ is bounded. 

\medskip {\sc Case 3} - We assume that $\delta_\eps\to 0$ and $\frac{\delta_{\epsilon}}{\sqrt{\epsilon}} \to +\infty$ as $\eps\to 0$. 

\medskip We let 
$$\check{\phi}_\eps(x)= \tilde{\phi}_{\epsilon} \left(\delta_\eps x+\tilde{x}_\eps \right) $$
for $x\in B_0\left(\delta \delta_\eps^{-1}\right)$. Then
$$\Delta_\xi \check{\phi}_\eps = \delta_\eps^2 \lambda_\eps e^{2\tilde{u}_\eps\left(\delta_\eps x+\tilde{x}_\eps\right)} \check{\phi}_\eps$$
in $B_0\left(\delta \delta_\eps^{-1}\right)$. Let $y_\eps\in \text{supp}\, \nu_{\eps}$ such that $d_g\left(x_\eps,y_\eps\right)=\delta_\eps$ and set 
$$\tilde{y}_\eps = \delta_\eps^{-1}\left( \exp_{g_l,x_l}^{-1}\left(y_\eps\right)-\tilde{x}_{\epsilon}\right)$$
so that 
\begin{equation}\label{eqcase3-1}
\tilde{y}_\eps\to \tilde{y}_0\hbox{ as }\eps\to 0
\end{equation}
after passing to a subsequence and set $R = \left|\tilde{y}_0\right|$. Thanks to the study of case 2, we know that 
\begin{equation}\label{eqcase3-2}
\check{\phi}_\eps\left(\tilde{y}_\eps\right) = O\left(1\right)\hskip.1cm.
\end{equation}
Thanks to the estimate (\ref{heat}) on the heat kernel, we also know that there exist $D_1>0$ and $r>0$ such that 
\begin{equation}\label{eqcase3-3}
\delta_\eps^2  e^{2\tilde{u}_\eps\left(\delta_\eps x+\tilde{x}_\eps\right)} \le D_1\hbox{ in }{\mathbb D}_{r}\left(0\right)\hskip.1cm.
\end{equation}
Assume first that $\check{\phi}_\eps$ does not vanish in ${\mathbb D}_{3R}\left(0\right)$. Then we can apply Harnack's inequality thanks to (\ref{eqcase3-3}) to get the existence of some $D_2>0$ such that 
\begin{equation}\label{eqcase3-4}
\check{\phi_\eps}(x)\ge D_2 \check{\phi}_\eps(0)
\end{equation}
for all $\eps>0$ and all $x\in {\mathbb D}_{\frac{r}{2}}(0)$. Note here that $\check{\phi}_\eps$ is maximal at $0$ thanks to the choice of $x_\eps$ we made. Since $\check{\phi}_\eps$ is super-harmonic on ${\mathbb D}_{\left|\tilde{y}_{\eps}\right|}\left(\tilde{y}_{\eps}\right)\subset {\mathbb D}_{3R}\left(0\right)$, we can also write that 
$$ \check{\phi}_\eps\left(\tilde{y}_\eps\right) \ge \frac{1}{2\pi \left\vert \tilde{y}_\eps\right\vert} \int_{\partial {\mathbb D}_{\left\vert \tilde{y}_\eps\right\vert}\left(\tilde{y}_\eps\right)} \check{\phi}_\eps\, d\sigma\hskip.1cm.$$
Keeping only the part of the integral which lies in ${\mathbb D}_{\frac{r}{2}}(0)$ and using (\ref{eqcase3-4}), we clearly get the existence of some $D_3>0$ such that 
$$\check{\phi}_\eps\left(\tilde{y}_\eps\right) \ge D_3 \check{\phi}_\eps(0)\hskip.1cm.$$
Here we used the assumption that $\check{\phi}_\eps>0$ in ${\mathbb D}_{3R}\left(0\right)$. Thanks to (\ref{eqcase3-2}), we conclude in this situation that $\check{\phi}_\eps(0)=\phi_\eps\left(x_\eps\right)=O(1)$. 

Assume now that $\check{\phi}_{\epsilon}$ vanishes in $\mathbb{D}_{3R}(0)$. By the claim \ref{nod}, since $\delta_{\epsilon} \to 0$ as $\epsilon \to 0$, $\check{\phi}_{\epsilon}$ also vanishes on $\partial \mathbb{D}_{4R}(0)$. By Cheng results on the nodal set of eigenfunctions (\cite{CHE1976}) and the Courant nodal theorem (\cite{COU1953}), the first eigenfunction $\check{\phi}_{\epsilon}$ vanishes on a piecewise smooth curve which connects two points $a_{\epsilon} \in \partial \mathbb{D}_{3R}(0)$ and $b_{\epsilon} \in \partial \mathbb{D}_{4R}(0)$. 

\medskip By \cite{ZIE1989}, corollary 4.5.3, there is a Poincar\'e inequality which bounds the $L^2$-norm of a function $\psi$ by the $L^2$ norm of its gradient with a multiplicative constant bounded by $B_{1,2}(\{\psi = 0\})^{-\frac{1}{2}}$ where $B_{1,2}$ is the Bessel capacity. The Bessel capacity is equivalent to the variational capacity (see \cite{TUR2000} Theorem 3.5.2) and we know that the variational capacity of a continuous curve which connects two uniformly distant points is uniformly bounded from below (see \cite{HEN2005bis}, pages 95-97).

\medskip Therefore, there is a constant $C$ independent of $\epsilon$ such that 
$$ \int_{\mathbb{D}_{4R}(0)}  \check{\phi}_{\epsilon}^2 dx \leq C \int_{\mathbb{D}_{4R}(0)} \left|\nabla\check{\phi}_{\epsilon}\right|^2 \, dx \hskip.1cm.$$
By the conformal invariance, the $L^2$-norm of the gradient is uniformly bounded. Thus $\left\{\check{\phi}_{\epsilon}\right\}$ is bounded in $L^2\left(\mathbb{D}_{\frac{r}{2}}(0)\right)$. By standard elliptic theory, thanks to (\ref{eqcase3-3}), we get also in this second situation that $\left\{\phi_{\epsilon}\left(x_{\epsilon}\right)\right\}$ is bounded. This ends the proof of the claim. \hfill $\diamondsuit$

\medskip We get now a quantitative non-concentration estimate on the $L^2$-norm of the gradient of $\phi_{\epsilon}^i$, $i=1,\dots,k$.

\begin{cl}\label{qgr}
There exists $C_2>0$ such that
$$\int_{B_g(x,r)} \left|\nabla \Phi_{\epsilon}\right|^2_g dv_g \leq \frac{C_2}{\sqrt{\log{\frac{1}{r}}}} $$
for all $\eps>0$ and all $r>0$. Here
$$ \left|\nabla \Phi_{\epsilon}\right|^2_g = \sum_{i=1}^k \left\vert \nabla \phi_\eps^i\right\vert_g^2\hskip.1cm.$$
\end{cl}

\medskip {\bf Proof.} It is clearly sufficient to prove the result for any $r$ small enough and any $x\in \omega_l$, $l$ fixed. Thus, setting as above
$$\tilde{\phi}_\eps^i = \phi_\eps^i\left(\exp_{g_l,x_l}(x)\right)$$
and ${\displaystyle \tilde{\Phi}_\eps =\left(\tilde{\phi}_\eps^i\right)_{i=1,\dots,k}}$, we need to prove that, for $r\le \delta$ and $x\in {\mathbb D}_\delta(0)$,
$$\int_{{\mathbb D}_r(x)}\left\vert \nabla  \tilde{\Phi}_\eps\right\vert_\xi^2\, dv_\xi \le \frac{C_2}{\sqrt{\ln\frac{1}{r}}}$$
for some $C_2>0$. In the following, we shall assume without loss of generality that $\delta<1$. Let us set 
$$F_\eps\left(r\right) =  \int_{{\mathbb D}_r(x)}\left\vert \nabla  \tilde{\Phi}_\eps\right\vert_\xi^2\, dv_\xi\hskip.1cm.$$
Using the equation satisfied by $\tilde{\Phi}_\eps$, namely
$$\Delta_\xi \tilde{\Phi}_\eps = \lambda_\eps e^{2u_\eps^l}\tilde{\Phi}_\eps$$
where 
$$u_\eps^l = \left(u_\eps+ v_l\right)\left(\exp_{x_l,g_l}\left(x\right)\right)\hskip.1cm,$$
we get that 
$$F_\eps\left(r\right)= \lambda_\eps \int_{{\mathbb D}_r(x)} e^{2u_\eps^l}\left\vert \tilde{\Phi}_\eps\right\vert^2\, dv_\xi + \int_{\partial {\mathbb D}_r(x)} \tilde{\Phi}_\eps\cdot \partial_\nu \tilde{\Phi}_\eps\, d\sigma_\xi\hskip.1cm.$$
Using now claims \ref{cap} and \ref{bou}, we can write that 
\begin{eqnarray*}
F_\eps\left(r\right)^2 &\le& \frac{D_1}{\left(\ln \frac{1}{r}\right)^2} + D_2 \left(\int_{\partial {\mathbb D}_r(x)} \left\vert \nabla \tilde{\Phi}_\eps\right\vert_\xi\, d\sigma_\xi\right)^2\\
&\le & \frac{D_1}{\left(\ln \frac{1}{r}\right)^2} + 2\pi r D_2 \int_{\partial {\mathbb D}_r(x)} \left\vert \nabla \tilde{\Phi}_\eps\right\vert_\xi^2\, d\sigma_\xi\\
&=&\frac{D_1}{\left(\ln \frac{1}{r}\right)^2} + 2\pi r D_2 F_\eps'\left(r\right)\hskip.1cm.
\end{eqnarray*}
We can write then that 
\begin{eqnarray*}
\left(F_\eps\left(r\right)\sqrt{\ln\frac{1}{r}}\right)'(s) &=& F_\eps'(s)\sqrt{\ln\frac{1}{s}} -\frac{1}{2s\sqrt{\ln \frac{1}{s}}} F_\eps(s)\\
&\ge & \frac{F_\eps(s)^2 \sqrt{\ln\frac{1}{s}}}{2\pi s D_2} - \frac{D_1}{2\pi s D_2 \left(\ln\frac{1}{s}\right)^{\frac{3}{2}}}-\frac{1}{2s\sqrt{\ln \frac{1}{s}}} F_\eps(s)\\
&\ge & -\frac{D_3}{s\left(\ln\frac{1}{s}\right)^{\frac{3}{2}}}
\end{eqnarray*}
for some $D_3>0$ independent of $x$, $\eps$ and $s$. Integrating this inequality from $r$ to $\delta$ leads to 
\begin{eqnarray*}
F_\eps\left(r\right)\sqrt{\ln\frac{1}{r}} &\le& F_\eps\left(\delta\right)\sqrt{\ln\frac{1}{\delta}} + \int_r^\delta \frac{D_3}{s\left(\ln\frac{1}{s}\right)^{\frac{3}{2}}}\, ds\\
&\le & \lambda_\eps \sqrt{\ln\frac{1}{\delta}} +\frac{2 D_3}{\sqrt{\ln \frac{1}{\delta}}}\\
\end{eqnarray*}
which clearly ends the proof of the claim. \hfill $\diamondsuit$

\medskip Thanks to the previous claims \ref{bou} and \ref{qgr}, we can compare precisely $\Phi_{\epsilon}$ and $K_{\epsilon}[\Phi_{\epsilon}]$~:

\begin{cl} \label{est} There exists $\beta_{\epsilon}\to 0$ as $\epsilon\to 0$ such that 
\begin{equation} \label{ineqn} \forall x \in M,\, \left\vert\Phi_{\epsilon}(x) \right\vert^2 \geq 1 - \beta_{\epsilon} \end{equation}
and
\begin{equation} \label{eqsu} \forall x \in \text{supp}(\nu_{\epsilon}),\, \left\vert K_{\epsilon}\left[\left\vert\Phi_{\epsilon} \right\vert\right](x) - 1\right\vert \leq \beta_{\epsilon} \end{equation}
\end{cl}

\medskip {\bf Proof.} We first prove that there exists $\beta_{\epsilon}\to 0$ as $\eps\to 0$ such that
\begin{equation} \label{epc} \forall x,y\in M, d_g(x,y) \leq \frac{\sqrt{\epsilon}}{\beta_{\epsilon}} \Rightarrow \left|\Phi_{\epsilon}(x) - \Phi_{\epsilon}(y) \right| \leq \beta_{\epsilon} \hskip.1cm.\end{equation}
For that purpose, let us set
$$\gamma_\eps = \left\Vert \eps e^{2u_\eps}\right\Vert_\infty^{\frac{1}{3}}\hskip.1cm.$$
Using claim \ref{cap} and (\ref{heat}), it is easily seen that $\gamma_\eps\to 0$ as $\eps\to 0$. Indeed, for any $r>0$,
$$\eps e^{2u_\eps(x)} \le \left(4\pi +o(1)\right)\int_{B_g\left(x,r\right)} d\nu_\eps + o(1) = 4\pi \nu\left(B_g\left(x,r\right)\right)+o(1)\le \frac{4\pi C_1}{\ln\frac{1}{r}} + o(1)\hskip.1cm.$$
We also have that 
$$\frac{\gamma_\eps}{\sqrt{\eps}}\to +\infty\hbox{ as }\eps\to 0$$
since $\gamma_\eps\ge \eps^{\frac{1}{3}}$. Let now $\left(x_\eps,y_\eps\right)\in M^2$ with $d_g\left(x_\eps,y_\eps\right)\le \frac{\sqrt{\eps}}{\gamma_\eps}$. Up to the extraction of a subsequence, $x_\eps\in \omega_l$ for some $l$ fixed. Let us set as before 
$$\tilde{\Phi}_\eps = \Phi_\eps\left(\exp_{g_l,x_l}(x)\right)$$
which satisfies 
$$\Delta_\xi \tilde{\Phi}_\eps = \lambda_\eps e^{2u_\eps^l} \tilde{\Phi}_\eps $$
with 
$$u_\eps^l = \left(u_\eps+v_l\right)\left(\exp_{g_l,x_l}(x)\right)\hskip.1cm.$$
We set 
$$\hat{\Phi}_\eps (x) = \tilde{\Phi}_\eps\left(\tilde{x}_\eps+\frac{\sqrt{\eps}}{\gamma_\eps}x\right)$$
where $x_\eps = \exp_{g_l,x_l}\left(\tilde{x}_\eps\right)$. We let $\alpha_\eps$ be the mean value of $\hat{\Phi}_\eps$ in ${\mathbb D}_2(0)$. By classical Sobolev and Poincar\'e inequalities, we know that there exists $D>1$ such that 
\begin{eqnarray*}
\left\Vert \hat{\Phi}_\eps-\alpha_\eps\right\Vert_{L^\infty\left({\mathbb D}_2(0)\right)} &\le &D \left\Vert  \Delta \hat{\Phi}_\eps\right\Vert_{L^\infty\left({\mathbb D}_2(0)\right)} + D \left\Vert \nabla \hat{\Phi}_\eps\right\Vert_{L^2\left({\mathbb D}_2(0)\right)}^2\\
&\le & D\lambda_\eps \gamma_\eps + D\frac{C_2}{\ln \frac{2\sqrt{\eps}}{\gamma_\eps}}\\
\end{eqnarray*}
thanks to claim \ref{qgr}. Setting 
$$\beta_\eps =2\left( D\lambda_\eps \gamma_\eps + D\frac{C_2}{\ln \frac{2\sqrt{\eps}}{\gamma_\eps}}\right)\hskip.1cm,$$
we then get that 
$$\left\vert \Phi_\eps\left(x_\eps\right)-\Phi_\eps\left(y_\eps\right)\right\vert \le \beta_\eps\hskip.1cm,$$
which clearly proves (\ref{epc}). 

\medskip We prove now that for all sequence $\{f_{\epsilon}\}$ of uniformly bounded functions which satisfy
\begin{equation} \label{epcf} \forall x,y\in M, d_g(x,y) \leq \frac{\sqrt{\epsilon}}{\beta_{\epsilon}} \Rightarrow \left|f_{\epsilon}(x) - f_{\epsilon}(y) \right| \leq \beta_{\epsilon} \hskip.1cm,\end{equation}
then, up to increase $\beta_{\epsilon}$, we have that
\begin{equation} \label{eh} \forall x\in M, \left|f_{\epsilon}(x) - K_{\epsilon}[f_{\epsilon}](x)\right| \leq \beta_{\epsilon}\hskip.1cm. \end{equation}
Indeed, for $x\in M$,
\begin{eqnarray*}
 \left|f_{\epsilon} - K_{\epsilon}[f_{\epsilon}]\right|(x) &\leq& \int_{B_{g}(x,\frac{\sqrt{\epsilon}}{\beta_{\epsilon}})} \left|f_{\epsilon}(x)-f_{\epsilon}(y)\right|p_{\epsilon}(x,y) dv_g(y) \\
 &&\quad+ \int_{M\setminus B_{g}(x,\frac{\sqrt{\epsilon}}{\beta_{\epsilon}})} \left|f_{\epsilon}(x)-f_{\epsilon}(y)\right|p_{\epsilon}(x,y) dv_g(y)
 \end{eqnarray*}
and by the property (\ref{epcf}) and estimates on the heat kernel,
$$\left|f_{\epsilon}(x) - K_{\epsilon}[f_{\epsilon}](x)\right| \leq \beta_{\epsilon} + 2\left\|f_{\epsilon}\right\|_{\infty} \left( O(e^{-\frac{1}{\beta_{\epsilon}^2}}) + O(e^{-\frac{1}{\sqrt{\epsilon}}})\right) $$
Up to increase $\beta_{\epsilon}$, we get (\ref{eh}).

\medskip Up to increase $\beta_{\epsilon}$, we get (\ref{eh}) for $f_{\epsilon} = \left|\Phi_{\epsilon}\right|^2$, thanks to (\ref{epc}). Then, by claim \ref{eul}, we easily get (\ref{ineqn}). By claim \ref{eul}, we also have that 
\begin{equation} \label{eqsup} \forall x\in \text{supp}(\nu_{\epsilon}), \left|\left|\Phi_{\epsilon}(x) \right|^2 - 1\right| \leq \beta_{\epsilon} \end{equation}
Again, up to increase $\beta_{\epsilon}$, we get (\ref{eh}) for $f_{\epsilon} = \left|\Phi_{\epsilon}\right|$ thanks to (\ref{epc}). Then, by (\ref{eqsup}), we easily get (\ref{eqsu}). This ends the proof of the claim. \hfill $\diamondsuit$

\medskip Thanks to claim \ref{est}, we can define $\Psi_{\epsilon} = \frac{\Phi_{\epsilon}}{\left|\Phi_{\epsilon}\right|} \in \mathcal{C}^{\infty}(M,\mathbb{S}^{k-1})$. Then, thanks to claim \ref{soboun}, $\{\Psi_{\epsilon}\}_{\epsilon}$ is bounded in $W^{1,2}\left(M,\mathbb{S}^{k-1}\right)$.

\begin{cl} \label{uni} There exists $C_3>0$ such that 
$$\left\vert \Psi_\eps(x)-\Psi_\eps(y)\right\vert^2 \sqrt{\ln \frac{2\delta(M)}{d_g\left(x,y\right)}}\le C_3$$
for all $x,y\in M$ and all $\eps>0$ where $\delta(M)$ is the diameter of $M$. In particular, the sequence $\{\Psi_{\epsilon}\}_{\epsilon}$ is uniformly equicontinuous in $\mathcal{C}^{0}(M,\mathbb{S}^{k-1})$.
\end{cl}

\medskip {\bf Proof.} We first claim that there exists $D_1>0$ such that 
\begin{equation} \label{ze} 
\sup_{x\in M}\sup_{v\in \Psi_{\epsilon}(x)^{\perp} \cap \mathbb{S}^{k-1}} \frac{1}{Vol_g\left(B_{g}\left(x,r\right)\right)}\int_{B_{g}(x,r)} \left( \Phi_{\epsilon} . v \right)^2 \, dv_g \leq \frac{D_1}{\sqrt{\ln{\frac{1}{r}}}} \end{equation}
for all $r$ small enough and all $\eps>0$.

For $x\in M$ and $v\in  \Psi_{\epsilon}(x)^{\perp} \cap \mathbb{S}^{k-1}$, the eigenfunction $\Phi_{\epsilon} . v$ vanishes at $x$. Using claim \ref{nod}, we can argue as in the proof of claim \ref{bou} to get the existence of some $D_2>0$ such that 
$$ \frac{1}{Vol_g\left(B_{g}\left(x,r\right)\right)}\int_{B_{g}(x,r)} \left( \Phi_{\epsilon} . v \right)^2 \, dv_g \le D_2 \int_{B_{g}(x,r)} \left\vert \nabla  \left( \Phi_{\epsilon} . v_{\epsilon} \right)\right\vert_g^2\, dv_g$$
for all $r$ small enough. We deduce thanks to claim \ref{qgr} that 
$$ \frac{1}{Vol_g\left(B_{g}\left(x,r\right)\right)}\int_{B_{g}(x,r)} \left( \Phi_{\epsilon} . v \right)^2 \, dv_g \le \frac{D_2 C_2}{\sqrt{\ln\frac{1}{r}}}$$
for all $r$ small enough and (\ref{ze}) follows.

\medskip Assume now by contradiction that the conclusion of the claim is false, that is there exists $\eps_n\to 0$ as $n\to +\infty$, $x_n$ and $y_n$ in $M$ such that 
\begin{equation}\label{nequlog}
\left\vert \Psi_{\eps_n}\left(x_n\right)-\Psi_{\eps_n}\left(y_n\right)\right\vert^2 \sqrt{\ln \frac{1}{r_n}}\to +\infty \hbox{ as }n\to +\infty
\end{equation}
where $r_n = d_g\left(x_n,y_n\right)\to 0$ as $n\to +\infty$. Thanks to (\ref{ineqn}) of claim \ref{est}, up to the extraction of a subsequence, there exists a fixed $v\in {\mathbb S}^{k-1}$ such that 
$$\frac{1}{Vol_g\left(B_g\left(x_n,r_n\right)\right)}\int_{B_g\left(x_n,r_n\right)} \left(\Phi_{\eps_n}\cdot v\right)^2\, dv_g \ge \frac{1-\beta_{\eps_n}}{k}=\frac{1}{k}+o(1)\hskip.1cm.$$
Thanks to (\ref{nequlog}), it is easy to find $X_n\in \Psi_{\eps_n}\left(x_n\right)^\perp$ and $Y_n\in \Psi_{\eps_n}\left(y_n\right)^\perp$ such that 
\begin{equation}\label{ze1}
v=X_n + Y_n \hbox{ and } \left\vert X_n\right\vert^2 + \left\vert Y_n\right\vert^2 = o\left(\sqrt{\ln\frac{1}{r_n}}\right)\hskip.1cm.
\end{equation}
We then write that 
\begin{eqnarray*}
\frac{1}{k}+o(1)&\le & \frac{1}{Vol_g\left(B_g\left(x_n,r_n\right)\right)}\int_{B_g\left(x_n,r_n\right)} \left(\Phi_{\eps_n}\cdot v\right)^2\, dv_g\\
&\le & \frac{2}{Vol_g\left(B_g\left(x_n,r_n\right)\right)}\int_{B_g\left(x_n,r_n\right)} \left(\Phi_{\eps_n}\cdot X_n\right)^2\, dv_g \\
&&+  \frac{2Vol_g\left(B_g\left(y_n,2r_n\right)\right)}{Vol_g\left(B_g\left(x_n,r_n\right)\right)}\frac{1}{Vol_g\left(B_g\left(y_n,2r_n\right)\right)}\int_{B_g\left(y_n,2r_n\right)} \left(\Phi_{\eps_n}\cdot Y_n\right)^2\, dv_g \\
&\le & 2 D_1 \left\vert X_n\right\vert^2 \left(\ln\frac{1}{r_n}\right)^{-\frac{1}{2}} + 16 D_1 \left\vert Y_n\right\vert^2 \left(\ln\frac{1}{2 r_n}\right)^{-\frac{1}{2}} \\
&=&o(1)\\
\end{eqnarray*}
using (\ref{ze}) and (\ref{ze1}). This is clearly a contradiction and proves the claim. \hfill $\diamondsuit$

\medskip Up to the extraction of a subsequence, one gets functions $\Phi \in W^{1,2}(M,\mathbb{R}^{k})\cap L^{\infty}(M,\mathbb{R}^{k})$ and $\Psi \in W^{1,2}(M,\mathbb{S}^{k-1})\cap \mathcal{C}^{0}(M,\mathbb{S}^{k-1})$ such that
\begin{equation}\label{phic} 
\Phi_{\epsilon} \rightharpoonup \Phi \hbox{ in } W^{1,2}(M,\mathbb{R}^k) \hbox{ and  }\Phi_{\epsilon} \rightarrow \Phi \hbox{ in } L^{p}(M,\mathbb{R}^k) \hbox{ as }\eps\to 0\end{equation}
and
\begin{equation}\label{psic} 
\Psi_{\epsilon} \rightharpoonup \Psi \hbox{ in } W^{1,2}(M,\mathbb{S}^{k-1}) \hbox{ and } \Psi_{\epsilon} \rightarrow \Psi \hbox{ in } \mathcal{C}^{0}(M,\mathbb{S}^{k-1})\hbox{ as }\eps\to 0 \end{equation}
where $\Psi$ and $\Phi$ satisfy
$$ \left|\Phi\right|^2 \geq_{a.e.} 1 \hbox{ and } \Psi = \frac{\Phi}{\left|\Phi\right|}\hskip.1cm. $$

\begin{cl} For $i\in \{1,\cdots,k\}$,
\begin{equation} \label{wlim} \phi^i_{\epsilon} e^{2u_{\epsilon}} dv_g \rightharpoonup^{*} \psi^i d\nu \hskip.1cm.\end{equation}
And, in a weak sense, we have that
\begin{equation} \label{limeq} \Delta_g \phi^i = \Lambda_1\left(M,[g]\right) \psi^i d\nu\hskip.1cm. \end{equation}
\end{cl}

\medskip {\bf Proof.} Let $\zeta \in \mathcal{C}^0(M)$. Then
\begin{eqnarray*}
 \int_M \zeta \phi_{\epsilon}^i e^{2u_{\epsilon}} dv_g - \int_M \zeta \psi^i d\nu &=&
 \int_M \left( K_{\epsilon}[\zeta\phi_{\epsilon}^i] - \zeta K_{\epsilon}[\phi_{\epsilon}^i] \right)d\nu_{\epsilon} \\
 &&+ \int_M \zeta \left( K_{\epsilon}[\phi_{\epsilon}^i] - \psi_{\epsilon}^i K_{\epsilon}[\left|\Phi_{\epsilon}\right|] \right)d\nu_{\epsilon}\\
 &&+\int_M \zeta \left(\psi_{\epsilon}^i K_{\epsilon}[\left|\Phi_{\epsilon}\right|] - \psi_{\epsilon}^i  \right)d\nu_{\epsilon}\\
 &&+\int_M \zeta \left( \psi_{\epsilon}^i d\nu_{\epsilon} - \psi^i d\nu \right)\hskip.1cm.
 \end{eqnarray*}
The first term converges to $0$ since $\{\phi_{\epsilon}^i\}$ is uniformly bounded thanks to claim \ref{bou}. The second term converges to $0$ since $\left(\left\vert \Phi_\eps\right\vert\right)$ is uniformly bounded thanks to claim \ref{bou} and $\{\psi_{\epsilon}^i\}$ is uniformly equicontinuous thanks to claim \ref{uni}. The third term converges to $0$ thanks to (\ref{eqsu}) (see claim \ref{est}). The last term also converges to $0$ thanks to the $C^0$-convergence of $\psi_{\epsilon}^i$ to $\psi^i$ (see (\ref{psic})) and the weak$^*$-convergence of $d\nu_\eps$ to $d\nu$. The first part of the claim follows. The second part of the claim is obtained by passing to the weak limit in the equations satisfied by the eigenfunctions thanks to (\ref{limlambdaepsilon}), (\ref{phic}) and (\ref{wlim}). \hfill $\diamondsuit$

\medskip We are now in position to end the proof of theorem \ref{thm1}. We test the equation (\ref{limeq}) against $\psi^i$ and sum over $i$ to obtain that
$$ \sum_{i=1}^k \int_{M} \left\langle \nabla\psi^i , \nabla\phi^i \right\rangle_g dv_g = \Lambda_1\left(M,[g]\right)\sum_{i=1}^k \int_M \left(\psi^i\right)^2 d\nu = \Lambda_1\left(M,[g]\right) \hskip.1cm.$$
Since
$$\nabla \psi^i= \nabla\left(\frac{\phi^i}{\left|\Phi\right|} \right) = \frac{\nabla \phi^i}{\left|\Phi\right|} - \frac{ \phi^i \nabla \left|\Phi\right|}{\left|\Phi\right|^2} \hskip.1cm,$$
we deduce that 
$$ \Lambda_1\left(M,[g]\right)=\sum_{i=1}^k \int_M\left\langle \nabla\psi^i , \nabla\phi^i \right\rangle_g dv_g = \int_M \left(\frac{\left|\nabla \Phi\right|_g^2}{\left|\Phi\right|} - \frac{\left|\nabla\left( \left|\Phi\right|\right)\right|_g^2}{\left|\Phi\right|}\right)dv_g\hskip.1cm.$$
Since $\Phi_{\epsilon} \rightharpoonup \Phi$ in $W^{1,2}(M,\mathbb{R}^k)$ and $\left|\Phi\right| \geq_{a.e.} 1$, we have the sequence of inequalities
\begin{eqnarray*} 
\Lambda_1\left(M,[g]\right) = \lim_{\epsilon \to 0} \int_M \left|\nabla \Phi_{\epsilon}\right|_g^2 dv_g &\geq& \int_M \left|\nabla \Phi\right|_g^2 dv_g \\
&\geq& \int_M \frac{\left|\nabla \Phi\right|_g^2}{\left|\Phi\right|}dv_g \\
&\geq& \Lambda_1\left(M,[g]\right) + \int_M \frac{\left|\nabla \left|\Phi\right|\right|_g^2}{\left|\Phi\right|}dv_g\\
&\ge &\Lambda_1\left(M,[g]\right) \hskip.1cm.
\end{eqnarray*}
Thus all the inequalities are in fact equalities and we deduce that $\left|\Phi\right|\equiv 1$ so that $\Psi=\Phi$ and that $\Phi_{\epsilon} \to \Phi$ in $W^{1,2}(M,\mathbb{R}^k)$ as $\eps\to 0$. We write that 
$$ 0 = \frac{1}{2}\Delta_g\left( \left|\Phi\right|^2 \right)= \sum_{i=1}^{k} \phi_i \Delta_g \phi_i - \sum_{i=1}^k \left|\nabla \phi_i\right|_g^2 = \Lambda_1\left(M,[g]\right) \left|\Phi\right|^2d\nu - \left|\nabla\Phi\right|_g^2$$
in a weak sense thanks to (\ref{limeq}) and what we just said. Then ${\displaystyle d\nu = \frac{ \left|\nabla\Phi\right|_g^2}{\Lambda_1\left(M,[g]\right)}dv_g} $ and the equation (\ref{limeq}) becomes
$$ \Delta_g \Phi = \left|\nabla\Phi\right|_g^2 \Phi $$
with $\Phi \in \mathcal{C}^0\left(M,\mathbb{S}^{k-1}\right) \cap W^{1,2}\left(M,\mathbb{S}^{k-1}\right)$. Then  $\Phi$ is weakly harmonic and by the regularity theory of harmonic functions by H\'elein (see \cite{HEL1996}), we can complete the proof of the theorem.

\section{Existence of maximal metrics for the first eigenvalue}\label{sectionminimalimmersion}

In this section, we prove theorem \ref{thm2}. Since it has already been proved in genus $0$ (Hersch \cite{HER1970}) and in genus 1 (Nadirashvili \cite{NAD1996}), we prove it for $\gamma\ge 2$. However, our proof clearly works in genus $1$ with light modifications (in the description of degeneracy of conformal classes) and this together with the result of El Soufi and Ilias \cite{ELS2000} give a new proof of the fact that the flat equilateral torus is maximizing the first eigenvalue of the Laplacian among the tori.

We let $M$ be a smooth compact orientable surface of genus $\gamma\ge 2$ and we let $\left(c_\alpha\right)$ be a sequence of conformal classes on $M$ such that 
\begin{equation}\label{eqmaxseq}
\lambda_\alpha = \Lambda_1\left(M,c_\alpha\right)\to \Lambda_1\left(\gamma\right)\hbox{ as }\alpha\to +\infty\hskip.1cm.\
\end{equation}
Let $h_\alpha$ be the hyperbolic metric of curvature $-1$ in the conformal class $c_\alpha$. By theorem \ref{thm1}, we know that there exists $g_\alpha\in c_\alpha$, smooth except at a finite set of conical singularities, such that 
\begin{equation}\label{eqmetrextr}
Vol_{g_\alpha}\left(M\right)=1\hbox{ and } \lambda_1\left(g_\alpha\right)=\Lambda_1\left(M,c_\alpha\right)=\lambda_\alpha\hskip.1cm.
\end{equation}
Moreover there exists a smooth harmonic map $\Phi_\alpha :\left(M,h_\alpha\right)\mapsto S^{k_\alpha}$ for some $k_\alpha\ge 2$ such that 
\begin{equation}\label{eqharmmap}
g_\alpha = \frac{\left\vert \nabla \Phi_\alpha\right\vert_{h_\alpha}^2}{\lambda_\alpha}h_\alpha\hskip.1cm.
\end{equation}
Since the multiplicity of eigenvalues is bounded by a constant which depends only on the genus $\gamma$ (see \cite{CHE1976}), the sequence $\left(k_\alpha\right)$ is uniformly bounded. Up to the extraction of a subsequence, we can assume in the following that $k_\alpha$ is fixed, $k_\alpha\equiv k$ for all $\alpha$.

\medskip The aim is to prove that the sequence $\left(h_\alpha\right)$ of hyperbolic metrics does converge smoothly to some hyperbolic metric as $\alpha\to +\infty$. In other words, we want to prove that the sequence of conformal classes $\left(c_\alpha\right)$ does not degenerate. For that purpose, we need to prove that the injectivity radius of $h_\alpha$ does not converge to $0$. Indeed, by Mumford's compactness theorem, a sequence of hyperbolic metrics with injectivity radius bounded from below does converge after passing to a subsequence (see \cite{MUM1971}) and the sequence of harmonic maps converges up to the formation of bubbles which correspond to points of concentration of the measure $dv_{g_\alpha}$ (see \cite{PAR1996}, \cite{SAC1981} or \cite{ZHU2010}, theorem 2.2). It is clear that claim \ref{claim-Kokarev} applies when the conformal class stays in a compact set so that such a concentration can not occur.  And theorem \ref{thm2} would follow.  

\medskip We proceed by contradiction and assume from now on that 
\begin{equation}\label{eqinjrad}
i_{h_\alpha}\left(M\right)\to 0\hbox{ as }\alpha\to +\infty\hskip.1cm.
\end{equation}
Then there exist $s$ closed geodesics $\gamma_\alpha^1,\dots,\gamma_\alpha^s$ whose length $l_\alpha^i$ goes to $0$ where $1\le s\le 3\gamma-3$ (see \cite{HUM1997}, IV, lemma 4.1). By the collar lemma (see \cite{HUM1997}, IV, proposition 4.2 or \cite{ZHU2010}, lemma 4.2, for the version we use), for all $1\le i\le s$, there exists an open neighborhood $P_\alpha^i$ of $\gamma_\alpha^i$ isometric to the following truncated hyperbolic cylinder
\begin{equation}\label{eqcylinder}
{\mathcal C}_\alpha^i = \left\{\left(t,\theta\right),\, -\mu_\alpha^i<t<\mu_\alpha^i,\, 0\le \theta<2\pi\right\}
\end{equation}
with 
\begin{equation}\label{eqdefmualphai}
\mu_\alpha^i = \frac{\pi}{l_\alpha^i}\left(\pi-2 \arctan\left(\sinh\frac{l_\alpha^i}{2}\right)\right)
\end{equation}
endowed with the metric
\begin{equation}\label{eqmetriccylinder}
h_\alpha^i = \left(\frac{l_\alpha^i}{2\pi \cos\left(\frac{l_\alpha^i}{2\pi}t\right)}\right)^2 \left(dt^2 +d\theta^2\right)\hskip.1cm.
\end{equation}
Note that we identify $\left\{\theta=0\right\}$ with $\left\{\theta=2\pi\right\}$ and that the closed geodesic $\gamma_\alpha^i$ corresponds to $\left\{t=0\right\}$. 

\medskip Let us denote by $M_\alpha^1,\dots,M_\alpha^r$ the connected components of ${\displaystyle M\setminus \bigcup_{i=1}^s P_\alpha^i}$. Then 
\begin{equation}\label{eqdecomposition}
M=\left(\bigcup_{i=1}^s P_\alpha^i \right)\bigcup \left(\bigcup_{j=1}^r M_\alpha^j\right)
\end{equation}
and this is a disjoint union. 

\medskip For $0<b<\mu_\alpha^i$, we let 
\begin{equation}\label{eqcyltrunc}
P_\alpha^i\left(b\right)=\left\{\left(t,\theta\right),\, -\mu_\alpha^i+b < t < \mu_\alpha^i -b\right\}
\end{equation}
after identification with ${\mathcal C}_\alpha^i$. We let also $M_\alpha^j\left(b\right)$ be the connected component of ${\displaystyle M\setminus \bigcup_{i=1}^s P_\alpha^i(b)}$ which contains $M_\alpha^j$. We claim that 

\begin{cl}\label{claimvolume}
There exists $D>0$ such that one of the two following assertions is true~:

\smallskip (a) There exists $i\in \left\{1,\dots,s\right\}$ such that 
$$Vol_{g_\alpha}\left(P_\alpha^i\left(a_\alpha\right)\right)\ge 1 -\frac{D}{a_\alpha}$$
for all sequences $a_\alpha\to +\infty$ with $\frac{a_\alpha}{\mu_\alpha^i}\to 0$ as $\alpha\to +\infty$ for all $1\le i\le s$.

\smallskip (b) There exists $j\in \left\{1,\dots,r\right\}$ such that 
$$Vol_{g_\alpha}\left(M_\alpha^j\left(9 a_\alpha\right)\right)\ge 1 -\frac{D}{a_\alpha}$$
for all sequences $a_\alpha\to +\infty$ with $\frac{a_\alpha}{\mu_\alpha^i}\to 0$ as $\alpha\to +\infty$ for all $1\le i\le s$.

\end{cl}

\medskip {\bf Proof.} We first construct test-functions for $\lambda_\alpha= \lambda_1\left(g_\alpha\right)$ compactly supported in the hyperbolic cylinders and in the $M_\alpha^j$'s. We let $b_\alpha\to +\infty$ as $\alpha\to +\infty$ with $\frac{b_\alpha}{\mu_\alpha^i}\to 0$ as $\alpha\to +\infty$ for all $1\le i\le s$.

\vspace{15mm}
\medskip {\underline{Test functions in the hyperbolic cylinders}.}

\medskip For $1\le i\le s$, we define $\varphi_\alpha^i$ as follows. It is $0$ outside of $P_\alpha^i$ and on $P_\alpha^i$, it is defined by 
\begin{equation}\label{eqdefphi}
\varphi_\alpha^i\left(t,\theta\right)=\left\{\begin{array}{cl}
{\displaystyle 0}&{\displaystyle \hbox{ for }-\mu_\alpha^i < t \le -\mu_\alpha^i+2b_\alpha}\\
\,&\, \\
{\displaystyle \frac{\mu_\alpha^i-2b_\alpha+ t}{b_\alpha}}&{\displaystyle \hbox{ for } -\mu_\alpha^i+2b_\alpha<t\le  -\mu_\alpha^i+3b_\alpha}\\
\,&\, \\
{\displaystyle 1}&{\displaystyle \hbox{ for }-\mu_\alpha^i+3b_\alpha<t<\mu_\alpha^i-3b_\alpha}\\
\,&\, \\
{\displaystyle \frac{\mu_\alpha^i-2b_\alpha- t}{b_\alpha}}&{\displaystyle \hbox{ for }\mu_\alpha^i-3b_\alpha\le t< \mu_\alpha^i-2b_\alpha}\\
\,&\, \\
{\displaystyle 0}&{\displaystyle \hbox{ for }\mu_\alpha^i-2b_\alpha\le t< \mu_\alpha^i}\\
\end{array}\right.
\end{equation}
We clearly have that 
\begin{equation}\label{eqtestfctgrad1}
\int_M \left\vert \nabla \varphi_\alpha^i\right\vert_{g_\alpha}^2\, dv_{g_\alpha} = \int_M \left\vert \nabla \varphi_\alpha^i\right\vert_{h_\alpha}^2\, dv_{h_\alpha}=\frac{4\pi}{b_\alpha}
\end{equation}
and that 
\begin{equation}\label{eqtestfctL21}
\int_M \left(\varphi_\alpha^i\right)^2\, dv_{g_\alpha}\ge Vol_{g_\alpha}\left(P_\alpha^i\left(3b_\alpha\right)\right)
\end{equation}
for $1\le i\le s$. 

\medskip {\underline{Test functions in the connected components $M_\alpha^j$}.}

\medskip For $1\le j\le r$, we define $\psi_\alpha^j$ as follows. It is $1$ in $M_\alpha^j$, $0$ in all the $M_\alpha^k$'s, $k\neq j$. And, in the $P_\alpha^i$'s, it is defined as follows. It is $0$ for $-\mu_\alpha^i+2b_\alpha\le t\le \mu_\alpha^i-2b_\alpha$. And then, for a given $i$, it depends~: if $\left\{t=\mu_\alpha^i\right\}$ is on the boundary of $M_\alpha^j$, then we let
$$\psi_\alpha^j =\left\{\begin{array}{cl}
{\displaystyle \frac{t+2b_\alpha-\mu_\alpha^i}{b_\alpha}}&{\hbox{ for }\mu_\alpha^i-2b_\alpha\le t \le\mu_\alpha^i-b_\alpha}\\
\,&\,\\
{\displaystyle 1}&{\hbox{ for }\mu_\alpha^i-b_\alpha\le t\le \mu_\alpha^i}\\
\end{array}\right.$$
Otherwise, we let $\psi_\alpha^j=0$ for $\mu_\alpha^i-2b_\alpha\le t \le \mu_\alpha^i$. We proceed in the same way to define $\psi_\alpha^j$ on the other side of the hyperbolic cylinder $P_\alpha^i$.

We clearly have that 
\begin{equation}\label{eqtestfctgrad2}
\int_M \left\vert \nabla \psi_\alpha^j\right\vert_{g_\alpha}^2\, dv_{g_\alpha} = \int_M \left\vert \nabla \psi_\alpha^j\right\vert_{h_\alpha}^2\, dv_{h_\alpha}=\frac{2\pi m_j}{b_\alpha}
\end{equation}
and that 
\begin{equation}\label{eqtestfctL22}
\int_M \left(\psi_\alpha^j\right)^2\, dv_{g_\alpha}\ge Vol_{g_\alpha}\left(M_\alpha^j\left(b_\alpha\right)\right)
\end{equation}
for $1\le j\le r$ where $m_j$ is the number of connected components of $\partial M_\alpha^j$. Note that $m_j\le 2\left(3\gamma-3\right)$.

\medskip For any two smooth functions $\varphi$ and $\psi$ on $M$ with disjoint compact supports, we have that 
$$\lambda_\alpha \le \max\left\{\frac{\int_M \left\vert \nabla\varphi\right\vert_{g_\alpha}^2\, dv_{g_\alpha}}{\int_M \varphi^2\, dv_{g_\alpha}}\, ;\,  \frac{\int_M \left\vert \nabla\psi\right\vert_{g_\alpha}^2\, dv_{g_\alpha}}{\int_M \psi^2\, dv_{g_\alpha}}\right\}\hskip.1cm.$$
Applying this to any pair of the above test functions, which all have disjopint compact supports, we get thanks to (\ref{eqtestfctgrad1}), (\ref{eqtestfctL21}), (\ref{eqtestfctgrad2}) and (\ref{eqtestfctL22}) that 
\begin{equation}\label{eqvol1}
\min\left\{Vol_{g_\alpha}\left(P_\alpha^i\left(3b_\alpha\right)\right)\,;\, Vol_{g_\alpha}\left(P_\alpha^j\left(3b_\alpha\right)\right)\right\}\le \frac{C}{b_\alpha}\hbox{ for }i\neq j\in \left\{1,\dots,s\right\}
\end{equation}
\begin{equation}\label{eqvol2}
\min\left\{Vol_{g_\alpha}\left(M_\alpha^i\left(b_\alpha\right)\right)\,;\, Vol_{g_\alpha}\left(M_\alpha^j\left(b_\alpha\right)\right)\right\}\le \frac{C}{b_\alpha}\hbox{ for }i\neq j\in \left\{1,\dots,r\right\}
\end{equation}
\begin{equation}\label{eqvol3}
\min\left\{Vol_{g_\alpha}\left(P_\alpha^i\left(3b_\alpha\right)\right)\,;\, Vol_{g_\alpha}\left(M_\alpha^j\left(b_\alpha\right)\right)\right\}\le \frac{C}{b_\alpha}\hbox{ for }1\le i\le s\hbox{ and }1\le j\le r
\end{equation}
where $C>0$ is some fixed constant independent of the sequence $\left(b_\alpha\right)$. 

\medskip Let $D>0$ that we shall fix later and let us assume that the conclusion of the claim does not hold. Let $\left(a_\alpha\right)$ be a sequence of positive real numbers with $a_\alpha\to +\infty$ and $\frac{a_\alpha}{\mu_\alpha^i}\to 0$ as $\alpha\to +\infty$ for all $1\le i\le s$. Assume by contradiction that for any $i\in\left\{1,\dots,s\right\}$,
\begin{equation}\label{eqvol4}
Vol_{g_\alpha}\left(P_\alpha^i\left(a_\alpha\right)\right)<1-\frac{D}{a_\alpha}
\end{equation}
and that, for any $j\in \left\{1,\dots,r\right\}$, 
\begin{equation}\label{eqvol5}
Vol_{g_\alpha}\left(M_\alpha^j\left(9a_\alpha\right)\right)<1-\frac{D}{a_\alpha}\hskip.1cm.
\end{equation}
Let $i\in \left\{1,\dots,s\right\}$. Assume that 
\begin{equation}\label{eqvol6}
Vol_{g_\alpha}\left(P_\alpha^i\left(3a_\alpha\right)\right)\ge \frac{10 C}{3a_\alpha}\hskip.1cm.
\end{equation}
Noting that $P_\alpha^i\left(3a_\alpha\right)\subset P_\alpha^j\left(a_\alpha\right)$, using (\ref{eqvol1}) with $3b_\alpha=a_\alpha$, we get that 
\begin{equation}\label{eqvol7}
Vol_{g_\alpha}\left(P_\alpha^j\left(a_\alpha\right)\right)\le \frac{3C}{a_\alpha}\hbox{ for }j\neq i\hskip.1cm.
\end{equation}
Using (\ref{eqvol3}) with $b_\alpha=a_\alpha$, we also get that 
\begin{equation}\label{eqvol8}
Vol_{g_\alpha}\left(M_\alpha^j\left(a_\alpha\right)\right)\le \frac{C}{a_\alpha}\hbox{ for }1\le j\le r\hskip.1cm.
\end{equation}
Since $Vol_{g_\alpha}\left(M\right)=1$, we deduce from (\ref{eqvol7}) and (\ref{eqvol8}) that 
$$Vol_{g_\alpha}\left(P_\alpha^i\left(a_\alpha\right)\right) \ge 1 - \frac{C\left(r+3s-3\right)}{a_\alpha}\hskip.1cm.$$
If we choose $D>C\left(r+3s-3\right)$, this contradicts (\ref{eqvol4}) and thus proves that (\ref{eqvol6}) can not hold. Thus, up to choose $D$ large enough, we have proved that 
\begin{equation}\label{eqvol9}
Vol_{g_\alpha}\left(P_\alpha^i\left(3a_\alpha\right)\right)\le \frac{10 C}{3a_\alpha} \hbox{ for }1\le i\le s\hskip.1cm.
\end{equation}

Let now $j\in \left\{1,\dots,r\right\}$ and assume that 
\begin{equation}\label{eqvol10}
Vol_{g_\alpha}\left(M_\alpha^j\left(3a_\alpha\right)\right)\ge \frac{2C}{3a_\alpha}\hskip.1cm.
\end{equation}
Since $M_\alpha^j\left(3a_\alpha\right)\subset M_\alpha^j\left(9a_\alpha\right)$, we can use (\ref{eqvol2}) with $b_\alpha= 9a_\alpha$ to write that 
\begin{equation}\label{eqvol11}
Vol_{g_\alpha}\left(M_\alpha^k\left(9a_\alpha\right)\right)\le \frac{C}{9a_\alpha}\hbox{ for }k\neq j\hskip.1cm.
\end{equation}
Using (\ref{eqvol3}) with $b_\alpha=3a_\alpha$, we can also write that 
\begin{equation}\label{eqvol12}
Vol_{g_\alpha}\left(P_\alpha^i\left(9a_\alpha\right)\right)\le \frac{C}{3a_\alpha}\hbox{ for }1\le i\le s\hskip.1cm.
\end{equation}
Combining (\ref{eqvol11}) and (\ref{eqvol12}) to the fact that $Vol_{g_\alpha}\left(M\right)=1$, we deduce that 
$$Vol_{g_\alpha}\left(M_\alpha^j\left(9a_\alpha\right)\right)\ge 1 - \frac{C\left(3s+r-1\right)}{9a_\alpha}\hskip.1cm.$$
Up to choose $D>\frac{C\left(3s+r-1\right)}{9}$, this contradicts (\ref{eqvol5}) and thus proves that (\ref{eqvol10}) can not hold. So we have proved that, up to choose $D$ large enough, 
\begin{equation}\label{eqvol13}
Vol_{g_\alpha}\left(M_\alpha^j\left(3a_\alpha\right)\right)\le \frac{2 C}{3a_\alpha} \hbox{ for }1\le j\le r\hskip.1cm.
\end{equation}
Now equations (\ref{eqvol9}) and (\ref{eqvol13}) together with the fact that $Vol_{g_\alpha}\left(M\right)=1$ and that $a_\alpha\to +\infty$ as $\alpha\to +\infty$ lead to a contradiction. Thus we have proved that (\ref{eqvol4}) and (\ref{eqvol5}) can not hold together, up to fix $D$ large enough. This clearly permits to end the proof of the claim. \hfill $\diamondsuit$

\medskip We shall now prove successively that both situations in claim \ref{claimvolume} lead to a contradiction.

\begin{cl}\label{claimcase1}
If (a) holds in claim \ref{claimvolume}, then $\Lambda_1\left(\gamma\right)\le 8\pi$. 
\end{cl}

\medskip {\bf Proof.} We follow ideas of Girouard \cite{GIR2009bis}. Let $a_\alpha\to +\infty$ with $\frac{a_\alpha}{\mu_\alpha^i}\to 0$ as $\alpha\to +\infty$ for all $1\le i\le s$. If (a) holds, there exists $1\le i\le s$ such that 
\begin{equation}\label{eqclaimcase1-1}
Vol_{g_\alpha}\left(P_\alpha^i\left(a_\alpha\right)\right) \ge 1-\frac{D}{a_\alpha}\hskip.1cm.
\end{equation}
Thus all the volume of $g_\alpha$ concentrates in the hyperbolic cylinder $P_\alpha^i$. We shall omit the subscript $i$ in the following and we shall identify $P_\alpha$ with ${\mathcal C}_\alpha$, a subset of $S^1\times {\mathbb R}$. We let $0\le \eta_\alpha\le 1$ be a smooth cut-off function defined on $M$ such that $\eta_\alpha\equiv 1$ on $P_\alpha^i\left(a_\alpha\right)$ and $\eta_\alpha \equiv 0$ on $M\setminus P_\alpha$. Moreover, we may choose it in such a way that 
$$\int_M \left\vert \nabla \eta_\alpha\right\vert_{g_\alpha}^2\, dv_{g_\alpha}\to 0\hbox{ as }\alpha\to +\infty$$
thanks to the fact that $a_\alpha\to +\infty$ as $\alpha\to +\infty$. We let $\Phi : P_\alpha\mapsto S^2$ be defined by 
$$\Phi\left(t,\theta\right) = \frac{1}{1+e^{2t}}\left(2e^t\cos \theta, 2e^t\sin\theta, e^{2t}-1\right)\hskip.1cm.$$
This map $\Phi$ is conformal. Thanks to Hersch (\cite{HER1970}, lemma 1.1), there exists a conformal diffeomorphism $\theta_\alpha$ of $S^2$ such that 
\begin{equation}\label{eqclaimcase1-2}
\int_{P_\alpha} \left(x\circ \theta_\alpha\circ \Phi\right)\eta_\alpha\, dv_{g_\alpha}=0\hskip.1cm.
\end{equation}
We let $i_\alpha\in \left\{1,2,3\right\}$ be such that 
\begin{equation}\label{eqclaimcase1-3}
\int_{P_\alpha\left(a_\alpha\right)} \left(x_{i_\alpha}\circ \theta_\alpha\circ \Phi\right)^2\eta_\alpha^2\, dv_{g_\alpha}\ge \frac{1}{3}\left(1-\frac{D}{a_\alpha}\right)
\end{equation}
and we set 
\begin{equation}\label{eqclaimcase1-4}
u_\alpha =\eta_\alpha\left( x_{i_\alpha}\circ\theta_\alpha\circ \Phi\right)\hskip.1cm.
\end{equation}
Such a $i_\alpha$ does obviously exist thanks to (\ref{eqclaimcase1-1}). It is then easily checked that 
$$\int_M \left\vert \nabla u_\alpha\right\vert_{g_\alpha}^2\, dv_{g_\alpha} \le \frac{8\pi}{3} +o(1)\hskip.1cm.$$
Then we have that 
$$\lambda_\alpha\le \frac{\int_M \left\vert \nabla u_\alpha\right\vert_{g_\alpha}^2\, dv_{g_\alpha}}{\int_M  u_\alpha^2\, dv_{g_\alpha}}\le 8\pi+o(1)\hskip.1cm.$$
This ends the proof of the claim.

\begin{cl}\label{claimcase2}
If (b) holds in claim \ref{claimvolume}, then $\Lambda_1\left(\gamma\right)\le \Lambda_1\left(\gamma-1\right)$. 
\end{cl}

\medskip {\bf Proof.} If (b) holds, there exists $1\le j\le r$ such that 
\begin{equation}\label{eqclaimcase2-1}
Vol_{g_\alpha}\left(M_\alpha^j\left(9a_\alpha\right)\right) \ge 1-\frac{D}{a_\alpha}\hskip.1cm.
\end{equation}
Thus all the volume of $g_\alpha$ concentrates in the connected component $M_\alpha^j$. We denote by $\tilde{M}_\alpha$ the connected component of $M\setminus \left(\gamma_\alpha^1\cup \dots\cup \gamma_\alpha^s\right)$ which contains $M_\alpha^j$. Then there exists a diffeomorphism $\tau_\alpha : \Sigma \mapsto \tilde{M}_\alpha$ with $\tau_\alpha^\star h_\alpha = \bar{h}_\alpha$ where $\left(\Sigma,\bar{h}_\alpha\right)$ is a hyperbolic surface (non-compact). We have that 
$$\bar{h}_\alpha \to h\hbox{ in }C^\infty_{loc}\left(\Sigma\right)\hbox{ as }\alpha\to +\infty\hskip.1cm.$$
We let for $\delta>0$
$$\Sigma_\delta = \left\{x\in \Sigma\hbox{ s.t. }i_h\left(x\right)\ge \delta\right\}$$
so that
$$h_\alpha\to h\hbox{ in }C^\infty\left(\Sigma_\delta\right)\hbox{ as }\alpha\to +\infty\hskip.1cm.$$
Up to a subsequence, there exist a decreasing sequence $\delta_\alpha\to 0$ and an increasing sequence $a_\alpha\to +\infty$ such that 
\begin{equation} \label{inclusion} M_\alpha^j\left(9a_\alpha\right) \subset \tau_\alpha\left(\Sigma_{\delta_\alpha}\right)\hskip.1cm. \end{equation}
We let $c=\left[h\right]$. We denote by $\left(\hat{\Sigma},\hat{c}\right)$ the compactification of the cusps of $\left(\Sigma,c\right)$ (see Hummel \cite{HUM1997} sections I.5, IV.2, IV.5 and V.1)~: $\left(\hat{\Sigma}\setminus \left\{p_1,\dots,p_t\right\},\hat{c}\right)$ is conformal to $\left(\Sigma,c\right)$. Note that $\hat{\Sigma}$ has genus less than or equal to $\gamma-1$. 

We also set 
$$\bar{\Phi}_\alpha = \Phi_\alpha\circ \tau_\alpha$$
and 
$$\bar{g}_\alpha = \tau_\alpha^\star g_\alpha\hskip.1cm.$$
We shall study the asymptotic behaviour of the harmonic map $\bar{\Phi}_\alpha : \left(\Sigma, \bar{h}_\alpha\right)\to {\mathbb S}^k$. By theorem 2.2 of Zhu \cite{ZHU2010}, there exist $x_1,\dots, x_N\in \Sigma$ and a harmonic map $\Phi : \left(\Sigma,h\right)\mapsto {\mathbb S}^k$ such that 
$$\bar{\Phi}_\alpha \to \Phi\hbox{ in }C^\infty_{loc}\left(\Sigma\setminus \left\{x_1,\dots, x_N\right\}\right)\hbox{ as }\alpha\to +\infty$$
and 
$$\int_{\tau_\alpha\left(\Sigma_{\delta_\alpha}\right)} \left\vert \nabla\Phi_\alpha\right\vert_{h_\alpha}^2\, dv_{h_\alpha}\to \int_\Sigma \left\vert \nabla \Phi\right\vert_h^2\, dv_h +\sum_{i=1}^N E_i$$
where the $E_i$'s correspond to the energies lost at the blow up points $x_i$. Since $\lambda_\alpha$ is uniformly bounded from below by $8\pi +\eps_0$ (because $\Lambda_1\left(\gamma\right)>8\pi$, see theorem \ref{thm1}), we can adapt claim \ref{claim-Kokarev} to prove that all the $E_i$'s are $0$. Now, thanks to theorem 3.6 of Sacks-Uhlenbeck \cite{SAC1981}, the harmonic map $\Phi$ can be extended to $\hat{\Sigma}$ by 
$$\hat{\Phi}:\left(\hat{\Sigma}, \hat{c}\right)\mapsto {\mathbb S}^k\hskip.1cm.$$
Choosing $g_0\in \hat{c}$ a regular metric, we have by conformal invariance of the $L^2$-norm of the gradient, (\ref{inclusion}), (\ref{eqclaimcase2-1}) and what we just said that 
$$\int_{\hat{\Sigma}} \left\vert \nabla \hat{\Phi}\right\vert_{g_0}^2\, dv_{g_0}=\Lambda_1\left(\gamma\right)\hskip.1cm.$$
We let 
$$g=\frac{\left\vert \nabla \hat{\Phi}\right\vert_{g_0}^2}{\Lambda_1\left(\gamma\right)} g_0$$
so that $Vol_g\left(\hat{\Sigma}\right)=1$. Let $\psi \in C^\infty\left(\hat{\Sigma}\right)$ be a first eigenfunction of $g$. Let also $\rho_\eps\in C^\infty_c\left(\hat{\Sigma}\setminus \left\{p_1,\dots,p_t\right\}\right)$ be such that 
$$\rho_\eps = 1 \hbox{ on }\hat{\Sigma}\setminus \bigcup_{i=1}^t B_{p_i}\left(\eps\right)$$
and such that 
$$\int_{\hat{\Sigma}}\left\vert \nabla \rho_\eps\right\vert_g^2\, dv_g\to 0\hbox{ as }\eps\to 0\hskip.1cm.$$
Then we write that 
$$\lambda_\alpha \le \frac{\int_{\Sigma} \left\vert \nabla\left(\rho_\eps\psi\right)\right\vert_{\bar{g}_\alpha}^2\, dv_{\bar{g}_\alpha}}{\int_\Sigma \left(\rho_\eps\psi\right)^2\, dv_{\bar{g}_\alpha}- \left(\int_{\Sigma}\rho_\eps\psi\, dv_{\bar{g}_\alpha}\right)^2}\hskip.1cm.$$
Passing to the limit as $\alpha\to +\infty$, we get that 
$$\Lambda_1\left(\gamma\right) \le \frac{\int_{\Sigma} \left\vert \nabla\left(\rho_\eps\psi\right)\right\vert_g\, dv_{g}}{\int_\Sigma \left(\rho_\eps\psi\right)^2\, dv_{g}- \left(\int_{\Sigma}\rho_\eps\psi\, dv_{g}\right)^2}\hskip.1cm.$$
Passing to the limit as $\eps\to 0$, it is easily checked that this leads to 
$$\Lambda_1\left(\gamma\right) \le \lambda_1\left(g\right)\le \Lambda_1\left(\gamma-1\right)\hskip.1cm.$$
This ends the proof of the claim.\hfill $\diamondsuit$

\medskip Thus we have proved that, if $\Lambda_1\left(\gamma\right)>\Lambda_1\left(\gamma-1\right)$, then $\Lambda_1\left(\gamma\right)$ is achieved by a smooth metric, up to a finite set of conical singularities. This ends the proof of theorem \ref{thm2}.

\section{The infimum of $\Lambda_1\left(M,[g]\right)$ over all conformal classes}\label{sectioninfsup}

\begin{figure}
\centering
\includegraphics[scale=0.6]{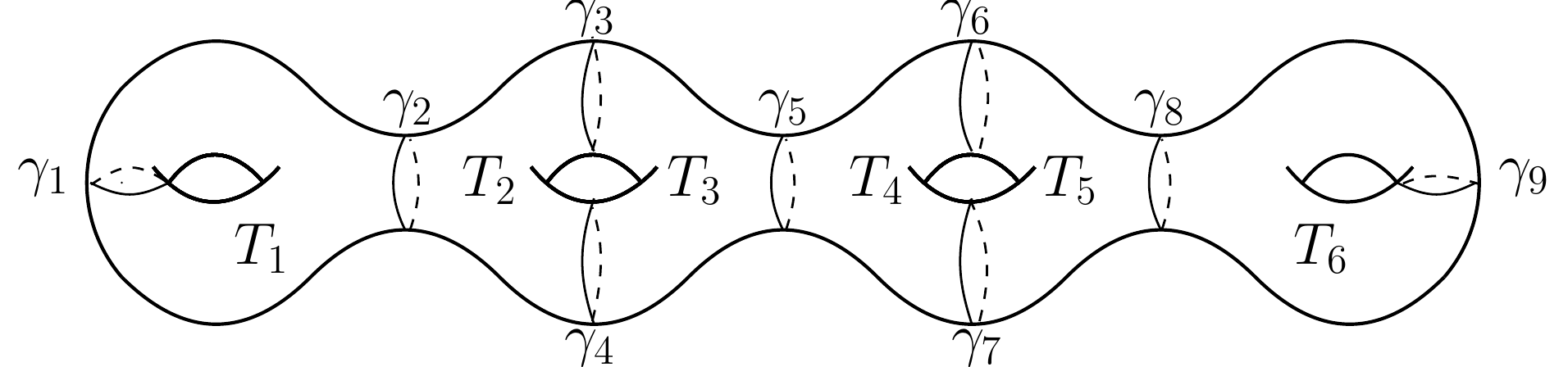}
\caption{Construction in genus $\gamma = 4$}
\label{const}
\end{figure}

In this section, we prove theorem \ref{thm3}. Fix $\gamma\ge 2$ since the result is already known in genuses 0 (\cite{HER1970}) and 1 (\cite{GIR2009bis}). We consider a sequence $M^n$ of hyperbolic surfaces (with metric $h_n$) obtained by gluing $2\gamma-2$ pairs of pants $T^j_n$~: these are surfaces containing $3\gamma-3$ closed geodesics $\gamma_n^1,\dots,\gamma_n^{3\gamma-3}$ of length $\eps_n\to 0$ as $n\to +\infty$ (see figure \ref{const}).

Each geodesic $\gamma_n^i$ has a neighbourhood $P_i^n$ isometric to the truncated cylinder $C_n = \left(-\nu_n,\nu_n\right)\times \left(0,2\pi\right)$ with 
$$\nu_n = \frac{\pi^2}{\eps_n} - \frac{2\pi}{\eps_n}\arctan\left(\sinh \frac{\eps_n}{2}\right)$$
endowed with the conformally flat metric
$$g_n = \left(\frac{\eps_n}{2\pi \cos\left(\frac{\eps_n}{2\pi}t\right)}\right)^2\left(dt^2+d\theta^2\right)\hskip.1cm.$$
We choose that the negative part of $P_i^n$, that is $t\le 0$, is in $T_k$ while the positive part is in $T_{k+1}$. 

We let $a_n\to +\infty$ with $\frac{a_n}{\nu_n} \to 0$ as $n\to +\infty$. We let $\psi$, $\varphi_l$, $\varphi_r$, $\theta_l$, $\theta_r$ be defined on $C_n$ as follows (depending on $n$ but we drop the subscript $n$)~:
\begin{eqnarray*}
\psi &=&\left\{\begin{array}{cl}
{\displaystyle 0}&{\displaystyle \hbox{ for }-\nu_n<t\le -\nu_n+a_n}\\
\,&\,\\
{\displaystyle \frac{t+\nu_n-a_n}{a_n}}&{\displaystyle \hbox{ for }-\nu_n+a_n<t\le -\nu_n+2a_n}\\
\,&\,\\
{\displaystyle 1}&{\displaystyle \hbox{ for }-\nu_n+2a_n<t<\nu_n-2a_n}\\
\,&\,\\
{\displaystyle \frac{\nu_n-a_n-t}{a_n}}&{\displaystyle \hbox{ for }\nu_n-2a_n\le t<\nu_n-a_n}\\
\,&\,\\
{\displaystyle 0}&{\displaystyle \hbox{ for }\nu_n-a_n\le t<\nu_n}\\
\end{array}\right.\\
&&\\
\varphi_l&=&\left\{\begin{array}{cl}
{\displaystyle 1}&{\displaystyle \hbox{ for }-\nu_n< t\le -\nu_n+a_n}\\
\,&\,\\
{\displaystyle \frac{2a_n-\nu_n-t}{a_n}}&{\displaystyle \hbox{ for }-\nu_n+a_n\le t<-\nu_n+2a_n}\\
\,&\,\\
{\displaystyle 0}&{\displaystyle \hbox{ for }-\nu_n+2a_n\le t<\nu_n}\\
\end{array}\right.\\
&&\\
\theta_l&=&\left\{\begin{array}{cl}
{\displaystyle \frac{t+\nu_n}{a_n}}&{\displaystyle \hbox{ for }-\nu_n\le t\le -\nu_n+a_n}\\
\,&\,\\
{\displaystyle 1}&{\displaystyle \hbox{ for }-\nu_n + a_n\le t\le -\nu_n+3a_n}\\
\,&\,\\
{\displaystyle \frac{4a_n-\nu_n-t}{a_n}}&{\displaystyle \hbox{ for }-\nu_n+3a_n\le t\le -\nu_n+4a_n}\\
\,&\,\\
{\displaystyle 0}&{\displaystyle \hbox{ for }-\nu_n+4a_n\le t\le \nu_n}\\
\end{array}\right.\\
\end{eqnarray*}
and $\theta_r(t)=\theta_l(-t)$, $\psi_r(t)=\psi_l(-t)$. 

\medskip We can now define the following test functions on $M^n$~: for $1\le i\le 3\gamma-3$, 
$$\psi_i = \left\{\begin{array}{ll}
{\displaystyle \psi}&{\displaystyle \hbox{ on }P_i}\\
\,&\,\\
{\displaystyle 0}&{\displaystyle \hbox{elsewhere}}\\
\end{array}\right.
$$
and 
\begin{eqnarray*}
\theta_{l,i} &=&\left\{\begin{array}{ll}
{\displaystyle \theta_l}&{\displaystyle\hbox{ on }P_i}\\
\,&\,\\
{\displaystyle 0}&{\displaystyle\hbox{ elsewhere}}\\
\end{array}\right.\\
&&\\
\theta_{r,i} &=&\left\{\begin{array}{ll}
{\displaystyle \theta_r}&{\displaystyle\hbox{ on }P_i}\\
\,&\,\\
{\displaystyle 0}&{\displaystyle\hbox{ elsewhere}}\\
\end{array}\right.\\
\end{eqnarray*}

For $1\le k\le \gamma-2$, 
$$\varphi_{2k+1} = \left\{\begin{array}{ll}
{\displaystyle \varphi_r}&{\displaystyle \hbox{ on }P_{3k}}\\
\,&\,\\
{\displaystyle \varphi_r}&{\displaystyle \hbox{ on }P_{3k+1}}\\
\,&\,\\
{\displaystyle \varphi_l}&{\displaystyle \hbox{ on }P_{3k+2}}\\
\,&\,\\
{\displaystyle 1}&{\displaystyle \hbox{ elsewhere in }T_{2k+1}}\\
\,&\,\\
{\displaystyle 0}&{\displaystyle \hbox{ elsewhere in }M^n}\\
\end{array}\right.
$$
and
$$\varphi_{2k} = \left\{\begin{array}{ll}
{\displaystyle \varphi_r}&{\displaystyle \hbox{ on }P_{3k-1}}\\
\,&\,\\
{\displaystyle \varphi_l}&{\displaystyle \hbox{ on }P_{3k}}\\
\,&\,\\
{\displaystyle \varphi_l}&{\displaystyle \hbox{ on }P_{3k+1}}\\
\,&\,\\
{\displaystyle 1}&{\displaystyle \hbox{ elsewhere in }T_{2k}}\\
\,&\,\\
{\displaystyle 0}&{\displaystyle \hbox{ elsewhere in }M^n}\\
\end{array}\right.
$$
We also define 
$$\varphi_1= \left\{\begin{array}{ll}
{\displaystyle \varphi_r+\varphi_l}&{\displaystyle \hbox{ on }P_1}\\
\,&\,\\
{\displaystyle \varphi_l}&{\displaystyle \hbox{ on }P_{2}}\\
\,&\,\\
{\displaystyle 1}&{\displaystyle \hbox{ elsewhere in }T_1}\\
\,&\,\\
{\displaystyle 0}&{\displaystyle \hbox{ elsewhere in }M^n}\\
\end{array}\right.
$$
and 
$$\varphi_{2\gamma-2}= \left\{\begin{array}{ll}
{\displaystyle \varphi_r+\varphi_l}&{\displaystyle \hbox{ on }P_{3\gamma-3}}\\
\,&\,\\
{\displaystyle \varphi_r}&{\displaystyle \hbox{ on }P_{3\gamma-2}}\\
\,&\,\\
{\displaystyle 1}&{\displaystyle \hbox{ elsewhere in }T_{2\gamma-2}}\\
\,&\,\\
{\displaystyle 0}&{\displaystyle \hbox{ elsewhere in }M^n}\\
\end{array}\right.
$$
We let now $g_n\in \left[h_n\right]$ with volume $1$ be such that 
$$\lambda_n = \lambda_1\left(g_n\right)=\sup_{g\in \left[h_n\right]} \lambda_1(g)Vol_g\left(M^n\right)\hskip.1cm.$$
Such a $g_n$ does exist thanks to theorem \ref{thm1}.

We denote by $\mathcal{E}$ the set of all the above functions defined on $M^n$. Note that all these functions have a $L^2$-norm (with respect to $g_n$) of their gradient converging to $0$ as $n\to +\infty$ (using the conformal invariance of this norm). Then, if $u$ and $v$ are two functions in $\mathcal{E}$ with disjoint compact supports, we have that 
\begin{equation} \label{disjointsupp} \lambda_n \min\left\{\int_{M^n} u^2\, dv_{g_n}\,;\,\int_{M^n} v^2\, dv_{g_n}\right\} \le o(1)\hbox{ as } n\to \infty \hskip.1cm. \end{equation}
Thanks to this remark, we will prove that one of the following situations must occur~:

\medskip a) Up to a subsequence, there exists $1\leq i \leq 3\gamma - 3$ such that 
$$\int_{M^n} \tau_i^2\, dv_{g_n}\to 1\hbox{ as }n\to +\infty$$
where 
$$\tau_i =\max\left\{\theta_{l,i},\theta_{r,i},\psi_i\right\}\hskip.1cm.$$

\medskip b) Up to a subsequence, there exists $1\le j\le 2\gamma-2$ such that 
$$\int_{M^n} \eta_j^2\, dv_{g_n}\to 1\hbox{ as }n\to +\infty$$
where for $1\leq k \leq \gamma-2$ we define :
$$ \eta_{2k+1} = \max\left\{  \varphi_{2k+1};\theta_{r,3k};\theta{r,3k+1};\theta_{l,3k+2}  \right\}$$
$$ \eta_{2k} = \max\left\{  \varphi_{2k};\theta_{r,3k-1};\theta{r,3k};\theta_{l,3k+1}  \right\}$$
$$ \eta_{1} = \max\left\{  \varphi_{1};\theta_{l,1};\theta{r,1};\theta_{l,2}  \right\} $$
$$ \eta_{2\gamma-2} = \max\left\{  \varphi_{2\gamma-2};\theta_{r,3\gamma-2};\theta_{l,3\gamma-3};\theta{r,3\gamma-3};  \right\} \hskip.1cm.$$

\medskip Indeed, we set
$$ \mathcal{F} = \{ u \in \mathcal{E} ; \int_M u^2 \not\rightarrow 0 \hbox{ as } n\to +\infty \} \hskip.1cm.$$
Since we have
$$ \int_{M^n} \left( \max_{u\in\mathcal{F}}\left\{ u^2 \right\} + \max_{u\in\mathcal{E}\setminus\mathcal{F}}\left\{ v^2 \right\} \right) dv_{g_n} \geq \int_{M^n} \max_{u\in\mathcal{E}}\left\{ u^2 \right\}dv_{g_n} = 1 $$
we easily get that
\begin{equation} \label{ineqconcentration} \int_{M^n} \left( \max_{u\in\mathcal{F}} u  \right)^2 dv_{g_n} \geq 1 - \sum_{v\in \mathcal{E}\setminus\mathcal{F}} \int_{M^n} v^2dv_{g_n} \rightarrow 1 \hbox{ as } n\to +\infty \hskip.1cm.\end{equation}
Then $\mathcal{F} \neq \emptyset$ and we distinguish two cases :

\medskip (i) There exists $1\leq j \leq 2\gamma-2$ such that $\varphi_j \in \mathcal{F}$.
Then, up to a subsequence, $\int_{M^n} \varphi_j^2 dv_{g_n}$ is uniformly bounded below and thanks to (\ref{disjointsupp}), we get that $\mathcal{F}$ contains at most two functions, with non-disjoint supports. Taking the maximum of these two functions, we easily obtain b) from (\ref{ineqconcentration}).

\medskip (ii) For all $1\leq j \leq 2\gamma-2$, $\varphi_j \in \mathcal{E} \setminus \mathcal{F}$. Since $\mathcal{F} \neq \emptyset$, there exists $1\leq i \leq 3\gamma - 3$ such that 
$$ \{ \psi_i ; \theta_{l,i} ; \theta_{r,i} \} \cap \mathcal{F} \neq \emptyset \hskip.1cm.$$
Then up to a subsequence, $\int_{M^n} \tau_i^2 dv_{g_n}$ is uniformly bounded below and thanks to (\ref{disjointsupp}), we get that $\mathcal{F} \subset \{ \psi_i ; \theta_{l,i} ; \theta_{r,i} \}$, and with (\ref{ineqconcentration}), we obtain a).

\medskip In both cases a) and b), we are in the situation of the lemma below and we deduce from it that $\lambda_n \le 8\pi+o(1)$. This concludes the proof of theorem \ref{thm3}. \hfill $\diamondsuit$

\medskip It remains to prove the following lemma we used during the previous proof~:

\begin{lem}\label{lemmaconcl}
Let $\Sigma$ be a compact orientable surface of genus $0$ with a boundary of $k$ connected components endowed with a sequence $g_n$ of metrics. Assume that there exists a sequence of functions $\eta_n : \Sigma\mapsto {\mathbb R}$ in $H^1\cap C^0$ such that~:

\smallskip i) $0\le \eta_n\le 1$.

\smallskip ii) $\eta_n$ is compactly supported in $\stackrel{\circ}{\Sigma}$.

\smallskip iii) ${\displaystyle \int_\Sigma \eta_n^2\, dv_{g_n}\to 1}$ as $n\to +\infty$.

\smallskip iv) ${\displaystyle \int_\Sigma \left\vert \nabla \eta_n\right\vert_{g_n}^2\, dv_{g_n}\to 0}$ as $n\to +\infty$.

\medskip Then there exists $u_n : \Sigma\mapsto {\mathbb R}$ in $H^1\cap C^0$ compactly supported in $\stackrel{\circ}{\Sigma}$ such that ${\displaystyle \int_\Sigma u_n\, dv_{g_n}=0}$ and 
$$\frac{\int_\Sigma \left\vert \nabla u_n\right\vert_{g_n}^2\, dv_{g_n}}{\int_\Sigma u_n^2\, dv_{g_n}}\le 8\pi +o(1)\hskip.1cm.$$
\end{lem}

\medskip {\bf Proof.}  
We first build a conformal diffeomorphism $\Psi_n : (\stackrel{\circ}{\Sigma},g_n) \rightarrow (\Sigma_n,h)$ where $\Sigma_n = \Psi_n(\stackrel{\circ}{\Sigma}) \subset \mathbb{S}^2$ and $h$ is the round metric of $\mathbb{S}^2$.

Let $U_1,\cdots,U_k$ some disjoint neighbourhoods of each connected component of the boundary which are diffeomorphic to annulus and such that, by the uniformization theorem for annuli (see \cite{HUM1997}, I.5), we get some conformal diffeomorphisms
$$ \Phi_n^i : (U_i,g_n) \rightarrow (A_{r_n^i},\xi) $$
where $0<r_n^i<1$ and for $0<r<1$, $A_r \subset \mathbb{D}$ is the annulus
$$ A_r = \{z\in\mathbb{C}; r<\left|z\right|<1\}\hskip.1cm.$$
Gluing $k$ copies of $\mathbb{D}$ instead of $A_{r_n^i}$, one can define a natural surface $\widetilde{\Sigma_n}$ endowed with a conformal structure $\widetilde{[g_n]}$ which extends $(\stackrel{\circ}{\Sigma},[g_n])$. $\widetilde{\Sigma_n}$ has a zero genus and by the uniformization theorem, there is a conformal diffeomorphism 
$$\widetilde{\Psi_n}: (\widetilde{\Sigma_n},\widetilde{[g_n]}) \rightarrow (\mathbb{S}^2,[h])\hskip.1cm.$$
Setting $\Psi_n = \widetilde{\Psi_n}_{|\stackrel{\circ}{\Sigma}}$ gives the expected conformal map. Following the arguments of claim \ref{claimcase1} permits to end the proof.
\hfill $\diamondsuit$

\bibliographystyle{plain}
\bibliography{biblio}

\end{document}